\documentclass[a4paper,12pt]{article}
\usepackage{amsmath}
\usepackage{amssymb, amsbsy, amstext}
\usepackage{srcltx}
\usepackage{amscd,amsfonts,color}
\newtheorem{lem}{Lemma}[section]%
\newtheorem{theorem}[lem]{Theorem}%
\newtheorem{prop}[lem]{Proposition}%
\newtheorem{conj}[lem]{Conjecture}%

\def\a{\alpha} \def\b{\beta} \def\g{\gamma} \def\d{\delta} \def\e{\varepsilon}
 \def\s{\sigma} \def\t{\tau}  
\def\th{\theta} \def\ld{\lambda} 
\def\Del{\Delta}  
\def\si{\Sigma} \def\O{\Omega} \def\G{\Gamma}
\def\D{{\rm D}}

\def\olg{\overline g}  \def\olm{\overline m}
 
\def\o{\overline}  
\def\di{\bigm|} \def\lg{\langle} \def\rg{\rangle}

\def\PSL{\hbox{\rm PSL}}  
\def\Aut{\hbox{\rm Aut\,}}  
 \def\soc{\hbox{\rm soc}} \def\Fix{\hbox{\rm Fix }}
  \def\mod{\hbox{\rm mod }}
  
 \def\PG{\hbox{\rm PG}} \def\PGL{\hbox{\rm PGL}}
  \def\GL{\hbox{\rm GL}}  \def\P\GL{\hbox{\rm P\GL}}
 \def\SL{\hbox{\rm SL}} \def\FF{{\hbox{\sf F\kern-.43emF}}}
  \def\Tr{\hbox{\rm Tr}}\def\Norm{\hbox{\rm Norm}}
\def\PGammaL{{\rm P\Gamma L}} \def\PSigmaL{{\rm P\Sigma L}}

\def\Sym{\hbox{\rm Sym}}
\def\Gal{\hbox{\rm Gal}}
\def\N{\hbox{\rm N}}
\def\C{\hbox{\rm C}}
\def\Z{\hbox{\rm Z}}
\def\A{\hbox{\rm A}}
\def\soc{\hbox{\rm soc}}
\def\u{\hbox{\rm u}}
\def\o{\hbox{\rm o}}
\def\J{\hbox{\rm J}}

\def\ZZ{\mathbb{Z}} 

\def\nd{\mathrel{\bigm|\kern-.7em/}} 
 \def\f{\noindent}
\def\qed{\hfill $\Box$} \def\demo{\f {\bf Proof}\hskip10pt}

\begin{document}
 \begin{center} {\bf\large On the Burness-Giudici Conjecture}
\vskip 3mm
{\sc Huye Chen}\\
{\footnotesize
China Three Gorges University,   Yichang 443002, P.R.China}\\

{\sc Shaofei Du}\\
{\footnotesize
Capital Normal University, Beijing 100048, P.R.China}\\

\end{center}

\renewcommand{\thefootnote}{\empty}%{footnote}}
  \footnotetext{E-mail addresses:  1120140003@mail.nankai.edu.cn, dushf@mail.cnu.edu.cn (S.F. Du).}

\begin{abstract}
Let $G$ be a permutation group on a set $\O$. A subset of $\O$ is a base  for $G$ if its pointwise stabilizer in $G$ is trivial. By  $b(G)$ we denote   the size of the smallest base of $G$. Every permutation  group with $b(G)=2$ contains some regular suborbits.
 It is conjectured   by  Burness-Giudici  in~\cite{BG} that  every primitive permutation group $G$ with $b(G)=2$ has the property  that  if $\a^g\not\in \G$ then
   $\G \cap \G^g\ne \emptyset$, where $\G$ is the union of all regular suborbits of $G$ relative to $\a$.
An affirmative answer  of the conjecture has been shown  for  many sporadic simple groups  and some  alternating groups   in \cite{BG}, but it  is still open for   simple groups of Lie-type. The first candidate of an infinite family of simple groups of Lie-type we should work on might be $\PSL(2,q)$, where $q\ge 5$. In this manuscript,  we show the correctness of the conjecture for all the primitive  groups with socle $\PSL(2,q)$, see Theorem~\ref{main}.
 \end{abstract}

\section{Introduction}
 A {\em base} for a finite permutation group $G$ on a set $\Omega$ is a subset $\Delta\subset \Omega$ whose pointwise  stabilizer is trivial.  The definition of a  base for finite permutation groups is a natural generalization of basis for  vector spaces. The {\em base size} of $G$ on $\O$, denoted by $b(G)$, is the minimal cardinality of a base for $G$. A {\em base size set} $\Delta$ for a finite permutation group $G\leq \Sym(\Omega)$,  is a base for $G$ on the set $\O$ with $|\Delta|=b(G)$.

 The first observation is that  a transitive group with $b(G)=1$ is just a regular group. So we are interested in the groups $G$ with $b(G)\ge 2$.

 The study of bases for finite permutation groups mainly includes two aspects: determining the base size for a finite permutation group and determining the base size set for a finite permutation group. In 1992, Blaha~\cite{Blaha} proved that finding the base size for a permutation group $G$ is NP-hard. After that, there are many researchers  working on three conjectures on $b(G)$, raised by Cameron, Kantor \cite{Cameron-Kantor}, Babai \cite{Babai-Cameron-Palfy} and Pyber \cite{Pyber}. Among them,  Pyber's conjecture was solved in 2018~\cite{DHM}.  All these conjectures  deal with the bounds of $b(G)$ for certain   primitive permutation groups,  while some  researchers also   give    descriptions of base size sets of such groups, see \cite{BGS1,BGS2,BBW,J1,J2,MNW}, for example.

Bases for a permutation group play an important role in the development of permutation group theoretic algorithms \cite{Seress}. Many of them are  related to   combinatorial structure, for an excellent survey article we refer to \cite{Bailey-Cameron}. Recently, Burness and Giudici \cite{BG} introduced a graph, called the {\em Saxl graph} which is related  to  bases of a permutation group with $b(G)=2$. Let $G\leq \Sym(\Omega)$ be a permutation group on $\Omega$ with $b(G)=2$. The vertex set of a Saxl graph $\Sigma(G)$ (denoted by $\Sigma$, simply) of $G$ on $\Omega$ is just $\Omega$ and two vertices are adjacent if and only if they form  a base size set. Clearly, if $G$ is a transitive permutation group, then the Saxl graph $\Sigma(G)$ is a vertex-transitive graph. Moreover, if $G$ is primitive, then $\Sigma(G)$ is connected,  but the converse is not true.
  In \cite{BG}, the authors  discussed  the valency, connectivity, hamiltonicity and the independence number of $\Sigma(G)$, and   proposed  the following conjecture about  the Saxl graph of a primitive permutation group.

\begin{conj}{\rm \cite{BG}} \label{BG}
Let $G$ be a finite primitive permutation group with $b(G)=2$ and Saxl graph $\Sigma(G)$. Then either $G$ is a Frobenius group and $\Sigma(G)$ is complete, or the diameter of $\Sigma(G)$ is $2$.
\end{conj}

 From now on, we assume $G$ is transitive on $\O$. Fix a point $\alpha\in \Omega$, the orbits of $G_\a$ are called {\em suborbits} of $G$ related to $\a$,  where $\{\a\}$ is said  to be trivial.
  The orbits of $G$ on $\Omega \times \Omega$ are called {\em orbitals} and each orbital  $\Delta$ corresponds to an orbital digraph for $G$ with the vertex set $\Omega$ and the arc set $\Delta$.
   Each orbital $\Delta$ corresponds to a suborbit   $\Delta(\alpha)=\{\beta\in \Omega \mid (\alpha,\beta)\in \Delta\}$ and  this   correspondence  between suborbits and orbitals
     is one to one.

  Now $\{\alpha,\beta\}$ is a base for $G$ if and only if $G_{\alpha}$ acts regularly on the suborbit containing $\b$ for $b(G)=2$. Therefore, the neighborhood $\si_1(\a)$ of $\alpha$ in $\Sigma(G)$ is the union $\G$  of the regular suborbits of $G$ relative to $\a$.
Let $d(\si)$ denote the diameter of the graph $\si$.  Note that  $G$ is a Frobenius group if and only if $d(\si)=1$. In this case,   $\O\setminus \{\a\}=\G$ and
     for any $\a^g\not\in \G$, we have $\a^g=\a$ and so $\G^g=\G$, or $\G\cap\G^g\ne\emptyset$.
   If $d(\si)=2$,  then for any $\a^g\not\in \G$,
       since the neighborhood of $\a^g$ is $\G^g$,   we have either $\a^g=\a$ and $\G=\G^g$ or $\a^g\ne \a$  and there exists a point $\b$ in  $\G\cap \G^g$, that is  $\G\cap\G^g\ne\emptyset$.
Conversely, suppose that for any    $\a^g\not\in \G$, we have   $\G\cap\G^g\ne\emptyset$. Then for any vertex  $\b=\a^g\not\in \G$, we have $\si_1(\a)\cap \si_1(\b)=\G\cap\G^g\ne\emptyset $, so $d(\si)\le 2$. Therefore,
     Conjecture~\ref{BG} is equivalent to  the following Conjecture~\ref{BG1}.

  \begin{conj} \label{BG1}
 Every primitive permutation group $G$ with $b(G)=2$ has the property  that  if $\a^g\not\in \G$ then
   $\G\cap \G^g\ne \emptyset$, where $\G$ is the union of all regular suborbits of $G$ relative to  a point $\a$.
\end{conj}

 Using a probabilistic approach, Burness and Giudici in \cite{BG} prove  the conjecture for some families of almost simple group. For example, the conjecture holds when $G=S_n$ or $G=A_n$ (with $n>12$) and the point stabilizer of $G$ is a primitive group.
Applying Theorem $1.4$ in \cite{Fawcett2} and Theorem $5.0.2$ in \cite{Fawcett1}, they also show the conjecture for some diagonal type and twisted wreath products primitive groups with sufficiently large order. By the Magma database they verified the Conjecture~\ref{BG} for all primitive groups of degree at most 4095. Also, the conjecture is shown for many sporadic simple groups in \cite{BG}.

So far this conjecture is still open for  simple groups of Lie-type. The first candidate of an infinite family of  these kinds of groups  we should work on might be $\PSL(2,q)$.  Therefore, the main goal of this manuscript is to show the correctness of Conjecture 1.2 for all the primitive permutation representations of the groups with socle $\PSL(2,q)$,  see  the following Theorem~\ref{main}, which will be proved in Section 3. 
\begin{theorem} \label{main} Let $q=p^n\ge  5$ for a prime $p$. Let $G$ be a primitive group with socle $\PSL(2,q)$ such that $b(G)=2$.  Then    $\G\cap \G^g\ne \emptyset$   if $\a^g\not\in \G$,
  where $\G$ is the union of all regular suborbits of $G$ relative to $\a$.
\end{theorem}
\section{Preliminary Results}
Let $G$ be a transitive permutation group on a finite set $\Omega$. For a subset $B\subset \Omega$, denote by $G_{B}$ and $G_{(B)}$ the subgroups of $G$ fixing $B$ set-wise and point-wise, respectively. If $B$ is a singleton $\{v\}$, then write $G_{v}=\{g\in G\mid v^g=v\}$, and call it the {\em stabilizer} of $v$ in $G$. For $v\in \Omega$, the orbit of $G$ containing $v$ is the subset $v^G:=\{v^g \mid g\in G\}$. Recall that $|v^G|=|G:G_{v}|$. If $G$ has only one orbit, then $G$ is said to be {\em transitive}. The permutation group $G$ is {\em semiregular} if $G_v=1$ for all $v\in \Omega$, and {\em regular} if further $G$ is transitive on $\Omega$.
 Let $H$ be a subgroup of $G$. By  $\N_G(H)$ and  $\C_G(H)$, we denote  the normalizer and  centralizer  of $H$ in $G$, respectively.
 By $\lceil a\rceil$, we denote the smallest integer no less than $a$.
 For any group $G$, set $\Z(G)$ denote the center of $G$. The results about the maximal subgroups of $G$ with socle $\PSL(2,q)$ were introduced in \cite{Giudici}.

\begin{prop} \label{psl} {\rm \cite[Chapter 3, Theorem 6.25]{Suz}} \ Let $q=p^n\ge  5$ be a power of  the prime $p$ and $d=(2,q-1)$. Then  a subgroup of $\PSL(2,q)$ is isomorphic to
\begin{enumerate}
\item[\rm(i)] $\PSL(2,p^m)$ where $m\di n$;
\item[\rm(ii)] $\PGL(2,p^m)$ where $2m\di n$ and $p$ is odd;
\item[\rm(iii)] $\D_{2(q\pm 1)/d}$ and their subgroups;
 \item[\rm(iv)] $\ZZ_p^n:\ZZ_{{\frac{q-1}d}}$ and their subgroups;
\item[\rm(v)] $A_5$, $q(q^2-1)\equiv 0(\mod 5)$;
\item[\rm(vi)] $S_4$, $q\equiv \pm 1(\mod 8)$;
  \item[\rm(vii)] $A_4$, either $p=2$ and $n$ is even  or $q$ is odd.
 \end{enumerate}
\end{prop}

\begin{prop} \label{max} {\rm \cite[Chapter 3, Ex.7, page 417]{Suz}} \ Let $q=p^n\ge  5$ be a power of  the prime $p$ and $d=(2,q-1)$.  Then every maximal subgroup of $\PSL(2,q)$ is one of the following:
\begin{enumerate}
\item[\rm(i)] $\D_{2(q+1)/d}$, where  $q\ne 7, 9$;
\item[\rm(ii)] $\PSL(2,p^m)$, where either  $\frac nm$ is an odd prime or $p=2$ and $n=2m$;
\item[\rm(iii)]  $\D_{2(q-1)/d}$, where $q\ne 5, 7, 9, 11$;
 \item[\rm(iv)] $A_5$, $q\equiv \pm 1(\mod 5)$ is a prime; or $q=p^2\equiv -1(\mod 5)$, where $p$ is an odd prime;
    \item[\rm(v)] $S_4$, $q\ge 5$ is an odd prime and $q\equiv \pm 1(\mod 8)$;
  \item[\rm(vi)] $A_4$, $q\ge 5$ is a prime and $q\equiv 3, 13, 27, 37(\mod 40)$;
  \item[\rm(vii)]  $\ZZ_p^n:\ZZ_{{\frac{q-1}d}}$;
   \item[\rm(viii)] $\PGL(2,p^m)$, where $n=2m$ and $p$ is an odd prime.
  \end{enumerate}
 \end{prop}

\begin{prop}{\rm \cite[p.407,(6.19)]{Suz}} \ \label{2,3} Let
$T=\PSL(2,q)$ and $G=\PGL(2,q)$. Set $d=(2,q-1)$.
\begin{enumerate}
\item[{\rm (1)}]   Let $Z$ be the center of the group $\SL (2,q).$
Then $x^2\in Z$ (resp. $x^3\in Z$) for $x\notin Z$ if and only if the trace
$\Tr(x)$ of $x$ is 0 (resp. $\pm 1$).
\item[{\rm (2)}] If $q$ is odd, then for any
involution $x\in T,$ $\C_T(x)\cong \D_{q\mp 1}$ and $\C_G(x)\cong
\D_{2(q\mp 1)},$ for $q\equiv \pm 1(\mod 4).$
\item[{\rm (3)}]
Acting on the projective line $\PG(1,q),$ every point-stabilizer
of $T$ (resp. $G$) is isomorphic to $\ZZ_p^n:Z_{\frac{q-1}d}$ (resp.
$\ZZ_p^n:Z_{q-1}$). Every element in a subgroup isomorphic to
$\ZZ_{\frac{q-1}d}$ (resp. $\ZZ_{q-1}$) of $T$ (resp.  $G$) fixes two
points.
\end{enumerate}
\end{prop}

%\begin{prop} \label{psl(2,q)} {\rm \cite[Chapter 3, Ex 3, page 416]{Suz}}  \  Let $j$ be the number of conjugacy classes of subgroups isomorphic to $\PSL(2,q)$ in $\PSL(2,q^n)$ where $q\ge 5$. Then $j=2$  when $q$ is odd and $n$ is even, while $j=1$ otherwise.
%\end{prop}

%\begin{prop} \label{A_4-S_4} {\rm \cite[Chapter 3, Ex 6, page 417]{Suz}}  If $q\equiv  \pm 1(\mod 8)$, then $\PSL(2,q)$ has two conjugacy classes of subgroups isomorphic to $A_4$ (resp. $S_4$).  If $q\equiv 3, 5(\mod 8)$,  then %all the subgroups isomorphic to $A_4$ are conjugate.
%\end{prop}

%\begin{prop} \label{A_5} {\rm \cite[Chapter 3, Ex 2, page 416]{Suz}}  \  The number of conjugacy classes of subgroups isomorphic to $A_5$ is 1,   if $q\equiv \pm 1(\mod 5)$ is even or if $q=5^{2m+1}$, and otherwise, the number is $2$.
%\end{prop}

\begin{prop}\label{man}
{\rm \cite{Man}} \
Let $G$ be a transitive group on $\O$ and let $H=G_\a$
for some $\a\in \O$.
Suppose that $K\le G$ and at least one $G$-conjugate of
$K$ is contained in $H$. Suppose further that the set of
$G$-conjugates of $K$ which are contained in $H$ form $t$
conjugacy classes of $H$ with representatives $K_1$, $K_2$, $\cdots,$ $K_t$.
Then $K$ fixes $\sum_{i=1}^{t}|\N_G(K_i):\N_H(K_i)|$ points of $\O$.
\end{prop}

The following $p$-group lemma  is  well-known and some related properties  are easy to check.
\begin{lem}  \label{maxcyclic}
Let $G$ be a finite nonabelian  $p$-group with order $p^n$ where $p$ is an odd prime and $n\ge 3$,  and one of whose  maximal  subgroups is cyclic. Then $G=\lg a,b\mid a^{p^{n-1}}=b^p=1, b^{-1}ab=a^{1+p^{n-2}} \rg$. Moreover, $G$ has only one  subgroup isomorphic to $\ZZ_p\times \ZZ_p$  and all the noncentral subgroups of order $p$ are conjugate to each other in $G$.
\end{lem}

Let $E=\FF_{q^n}$, $F=\FF_q$ and $E/F$ be Galois. Set $G=GalE/F=\{\eta_1=1,\eta_2,\cdots,\eta_n\}$ to be a Galois group. If $u\in E$, then define
$$Tr(u)=\sum_1^n\eta_i(u), ~~Norm(u)=\prod_1^n\eta_i(u)$$
and call these respectively the {\rm trace} and {\rm norm} of $u$ in $E/F$. Evidently, $Tr(\eta(x))=Tr(x)$ and $Norm(\eta(y))=Norm(y)$ for any $x,y\in E$.

%Let $E=\FF_q^n$ and $F=\FF_q$. For $x \in E$, we define the {\em trace} of  $x$ over $F$ as $Tr(x)=x+x^q+\cdots+x^{q^{n-1}}$. The {\em norm} over $F$ of an element $y\in E$ is defined by $Norm(y)=yy^q\cdots y^{q^{n-1}}=y^{\frac{q^n-1}{q-1}}$. It is easy to say that $Tr(x^q)=Tr(x)$ and $Norm(y^q)=Norm(y)$ for any $x,y\in E$. The following two results are useful during our proof.

\begin{prop} (Hilbert's Satz 90) \label{Gal}
{\rm \cite[Theorem 4.29,4.30,4.33]{Jac} }Let $E/F$ be cyclic with Galois group $G=\lg \eta\rg $.
  Let $d_1, d_2\in E$ such that the trace of $d_1$ is 0 and the norm of $d_2$ is 1.   Then there exist   $c_1, c_2\in E$  such that $d_1=c_1-\eta (c_1)$ and
  $d_2=c_2\eta(c_2)^{-1}$.
 \end{prop}

\begin{prop} (Hilbert's Satz 90) \label{Gal2}
{\rm \cite[Theorem 4.36]{Jac}} Let $F$ be a field of characteristic   $p\ne 0$ and  let $E/F$ be a $p$-dimensional cyclic extension with Galois group $G=\lg \eta\rg $. Then $E=F(c)$, where $\eta(c)=1+c$ and $c^p-c\in F$.
 \end{prop}

\vskip 3mm  The following proposition   is modified from   Lemma $3.8$ and Lemma $3.6$ in \cite{BG}.
\begin{prop}\label{smallp}
 Consider the transitive permutation representation of $G$ on $\O:=[G:M]$. Let $\hat{Q}(G):=\frac{|M|}{|\O|}\sum_{i=1}^k(|x_i|-1)\frac{|\Fix(\lg x_i \rg)|}{|\N_M(\lg x_i\rg )|}$, where $\mathcal{P}^*(M):=\{\lg x_1\rg, \lg x_2\rg,\cdots, \lg x_k\rg\}$ is the set of representatives  of conjugacy classes of subgroups of prime order in $M$ and $\Fix(\lg x_i\rg)$ is the set of fixed points of $\lg x_i\rg$ on $\O$.  If $\hat{Q}(G)<\frac{1}{2}$, then any two vertices in $\Sigma(G)$ have a common neighbor.
  \end{prop}

The following  result about Saxl graphs is clear.
\begin{lem}  \label{subgraph} Let $G$ be a transitive group on $\O$ and $b(G)=2$. Suppose that the Saxl graph $\si(G)$ has diameter $2$. Then for any nonregular transitive subgroup $G_1$ on $\O$,
we have $b(G_1)=2$ and the Saxl graph of $G_1$ has diameter $2$.
\end{lem}
\demo Let $\{ \a, \b\}$ be a base of $G$. Then $(G_1)_{(a,\b)}\le G_{(\a,\b)}=1$ and so $\{\a, \b\}$ is a base of $G_1$, that is $b(G_1)=2$.
Let $\G$ (resp. $\G_1$) be the respective   union of all regular suborbits of $G$ (resp. $G_1$)  relative to $\a$. Then  $\G$ is also a union of some regular suborbits of $G_1$ and so $\G\subset \G_1$. For any $\b\in \O\setminus \G$, since $G_1$ is transitive on $\O$, there exists
$g\in G_1\le G$ such that $\b=\a^g$. Now, $\emptyset \ne  \G\cap \G^g\subset \G_1\cap \G_1^g$, and the Saxl graph of $G_1$ has diameter $2$.\qed

%Section 3
\section{Proof of Theorem \ref{main}}
From now on, let $q=p^n\ge  5$ for a prime $p$ and an integer $n$, and  let  $\FF_q$  be a finite field of order $q$. Let $\s \in \Aut(\FF_q)$ defined by $\s(x)=x^p$,
 for any $x\in \FF_q$. Let $T=\PSL(2,q)$, a nonabelian simple group and  $T\le G\le \PGammaL(2,q)=(\PSL(2,q):\lg \delta \rg):\lg f\rg $, where $f$ is induced by $\s$ and $\lg \delta \rg\cong \ZZ_{(2,q-1)} $.

 To show Theorem~\ref{main},
  let $M$  be  any maximal subgroup  of $G$ and $\O=[G:M]$ which is the set of right cosets of $M$ in $G$. Then we consider the  right multiplication action of $G$ on $\O$, which is primitive. Then the Saxl graph  $\si$ of $G$  is connected.  As before, by $\G$ we denote the union of regular suborbits of $G$ relative to $\a=M$.
  To prove Theorem~\ref{main}, it suffices to   show the diameter $d(\si)$  of $\si$  is 2.
  What we need to do is to analyse the permutation group $G$ on $\O$ relative to  every  maximal subgroup $M$ up to conjugacy.
 Denoted by $M_G$ the kernel of this action. Suppose that $M_G\ne 1$.  Then $T\le M_G$, and  $G/M_G$ is abelian.  In this case, $b(G/M_G)=1$. Therefore, we assume $M_G=1$, that is $T\not\le M$.   Set $M_0=M\cap T$. Then it is proved in  \cite{Giudici} that $M_0$ is maximial in $T$, except for the groups $G$ in Figure~\ref{not maximal}.
 The first nine cases in this figure can be checked by Magma that their Saxl graphs have diameter $2$ for $b(G)=2$. The tenth case will be done in Lemma \ref{GLD}.
 Then all the cases for $q=7,9$ (resp. $q=7,9,11$) with $M_0=\D_{q+1}$ (resp. $M_0=\D_{q-1}$) have been checked by Magma. In fact, they are all in Figure \ref{not maximal}. 
 Moreover, we have $M/M_0\cong G/T$ by \cite{Giudici}. By $\Fix(K)$, we denote the set of fixed points of $G$ acting on $\O=[G:M]$. In what follows, the eight cases in Proposition~\ref{max}  will be treated   in the following  four subsections, separately.

  \begin{figure}

        \centering

\begin{tabular}{|c|c|}

        \hline

  $G$ & $M$  \\
\hline
\hline

         $\PGL(2,7)$& $\N_G(\D_6)=\D_{12}$ \\

         $\PGL(2,7)$ & $\N_G(\D_8)=\D_{16}$\\
        $\PGL(2,9)$& $\N_G(\D_{10})=\D_{20}$ \\

         $\PGL(2,9)$ & $\N_G(\D_8)=\D_{16}$\\
         $\lg \PSL(2,9),\delta f\rg$& $\N_G(\D_{10})=\ZZ_5:\ZZ_4$ \\

        $\lg \PSL(2,9),\delta f\rg$& $\N_G(\D_{8})=\ZZ_8:\ZZ_2$\\
        $\PGammaL(2,9)$& $\N_G(\D_{10})=\ZZ_{10}:\ZZ_4$ \\

         $\PGammaL(2,9)$ & $\N_G(\D_8)$\\

         $\PGL(2,11)$& $\N_G(\D_{10})=\D_{20}$ \\

         $\PGL(2,q)$, $q=p\equiv \pm 11,19(\mod 40)$ & $\N_G(A_4)=S_4$\\
         \hline
 \end{tabular}

       \caption{$M_0$ is not maximal in $\PSL(2,q)$ for $T\nleq M$}
       \label{not maximal}
\end{figure}

\subsection{$M_0=\D_{\frac{2(q+1)}{d}}$ where $d=(2,q-1)$}

\vskip 3mm
The main result of this subsection is the following Theorem~\ref{main1}. And we assume that $q\neq 7,9$ in this subsection.
\begin{theorem}\label{main1}
Suppose that $T=\soc(G)=\PSL(2,q)$ where $q\ge 5$ and $q\neq 7,9$.  Set $M$ to be a maximal subgroup of $G$ with $M_0=M\cap T=\D_{q+1}$. Consider the primitive permutation representation of $G$ on $\O=[G:M]$. Then the following statements hold.
\begin{enumerate}
\item[\rm(i)] We have $b(G)=2$ if and only if $q$ is odd and $G\in\{T:\lg f^j\rg, T.\lg \delta f^i\rg \}$ with $0\le j \le n-1$, $1\le i \le n-1$ and $\frac {n}{(n,i)}$ is even.
\item[\rm(ii)]  The Saxl graph $\si(G)$ has diameter $2$ provided $b(G)=2$.
\end{enumerate}
 \end{theorem}
In Subsection 3.1.1, the action of $\PSL(2,q)$ on its cosets relative to $M=\D_{q+1}$   will be characterized, where $q$ is odd;  in Subsection 3.1.2, it will be  shown that  $\si(G)$ has diameter 2 where $G=\PSigmaL(2,q)$ relative to $\D_{q+1}.\ZZ_n$ for odd $q$; and  finally in Subsection 3.1.3,  Theorem~\ref{main1} will be proved.

 \subsubsection{$G=\PSL(2,q)$ and $M=D_{q+1}$ for $q$ is odd}
Let $G=\PSL(2,q)$ and $M=\D_{q+1}$ where $q=p^n$ is odd.
Then  $G$ has only one conjugacy class of subgroups isomorphic to $\D_{q+1}$, while if $q\ne 7, 9$, they are maximal.  Pick up one of such subgroups $\a=M$ and set $\O=[G:M]$.

For a finite field $\FF_q$, set $S^*={\FF_q^*}^2$, $S=\{0\}\cup S^*$ and   $N=\FF_q\setminus S$.
 First,  we need a technical lemma proved by Prof. Keqin Feng.

\begin{lem}\label{feng}
Suppose that $q\ge 17$ is odd. Then for any $t\in \FF^*_q\setminus\{1\}$, set
$n=|(1+N)\cap (t+N)\cap S^*|$. Then
$$n\ge \lceil\frac{1}{8}(q-2\sqrt{q}-7)\rceil.$$
\end{lem}
\demo
Set $\eta: \FF^*_q\rightarrow\{\pm1\}$ by assigning the elements of $S^*$ to $1$ and that of $N$ to $-1$ and moreover, set $\eta(0)=0$. This $\eta$ is exactly that in \cite[Example 5.10]{LN}. Also we need to quote the following three results from~\cite[Theorem 5.4, 5.48, 5.41]{LN}:
\begin{enumerate}
\item[\rm(i)] $\sum_{x\in \FF_q}\eta(x)=0$;
\item[\rm(ii)] $\sum_{x\in \FF_q}\eta(x^2+Ax+B)=q-1$ for $A^2-4B=0$ or $-1$ for otherwise, where $A,B\in \FF_q$;
\item[\rm(iii)] Set $m:=\sum_{x\in \FF_q}\eta(x(x-1)(x-t))$, where $t\in \FF_q$. Then $|m|\leq 2\sqrt{q}$.
\end{enumerate}

Set $W=(1+N)\cap (t+N)\cap S^*$, that is $$W=\{x\in\FF_q\mid \eta(x-1)=\eta(x-t)=-1, \eta(x)=1\}.$$
Now let $t\ne 1$.
Then by the above three formulas, we have
$$\begin{array}{lcl}
 |W|&=&\frac{1}{8}\sum_{x\in \FF_q\setminus\{0,1,t\}}(1+\eta(x))(1-\eta(x-1))(1-\eta(x-t))\\
     &=&\frac{1}{8}\sum_{x\in \FF_q\setminus\{0,1,t\}}(1+\eta(x)-\eta(x-1)-\eta(x-t)-\eta(x(x-1))-\eta(x(x-t))\\
     &&+\eta((x-1)(x-t))+\eta((x-1)(x-t)x)\\
     &=&\frac{1}{8}[(q-3)+(-1-\eta(t))+(\eta(-1)+\eta(t-1))+(\eta(-t)+\eta(1-t))\\
     &&+(1+\eta(t(t-1)))+(1+\eta(1-t))+(-1-\eta(t))+m]\\
     &=&\frac{1}{8}[q-3+m+l],
\end{array}$$
where
$$l=\begin{cases}
3\eta(t-1)+\eta(t^2-t)+1-\eta(t),~if~ \eta(-1)=1 ,\\
\eta(1-t)+\eta(t^2-t)-3\eta(t)-1,~if~ \eta(-1)=-1 .
\end{cases}$$
Then $|l|\leq 4$ and
$$|W|\geq \frac{1}{8}[(q-3)-2\sqrt{q}-4]=\frac{1}{8}[(\sqrt{q}-1)^2-8].$$
So if $\sqrt{q}-1>\sqrt{8}$, that is $q\geq17$, then we have $W\neq \emptyset$.
\qed

\begin{lem}\label{connected}
Let $G=\PSL(2,q)$, $M\le G$ and $M\cong \D_{q+1}$ where $q\ge 17$ is odd. Then every  involution $\olg :=gZ\in G$ has the following form
 $g={\small \left(
\begin{array}{cc}
a&b\\
c&-a \\
\end{array}
\right)},$ where $Z$ is the center of the group $\SL (2,q)$ and
 \begin{eqnarray}
\label{abc}
a^2+bc=-1.
\end{eqnarray}
Moreover, for any $\overline{g}\in G\setminus M$,  there exist at least  $2\lceil\frac{1}{8}(q-2\sqrt{q}-7)\rceil$ noncentral  involutions  $\olm $ in $M$ such that $|\olm\olg|\di \frac{q+1}2$ and $|\olm \olg|>2$.
\end{lem}
\demo
Set $\FF_{q}^*=\lg \th \rg $ and let $\FF_{q^{2}}=\FF_{q}(\zeta)$ be the extension of $\FF_q$ by adding a  root $\zeta$ of the irreducible polynomial   $x^2-\th$ over $\FF_q$. Then every element in $\FF_{q^2}$ can be uniquely written as $u_0+v_0\zeta$ for some $u_0,v_0\in\FF_q$.  Moreover,  $\FF_{q^2}$ can be viewed as a $2$-dimensional vector space over $\FF_q$ and every element $u_0+v_0\zeta$ gives a linear transformation:
 $a\to (u_0+v_0\zeta) a, \, \forall a\in \FF_{q^2}$, which corresponds to the matric ${\small \left(
\begin{array}{cc}
u_0&v_0\theta \\
v_0&u_0\\
\end{array}
\right)}$   under the base: $1, \zeta$. This gives an embedding of $\FF_{q^2}^*$ into $\GL(2,q)$.

%Take a generator $u_0+v_0\zeta$ of $\FF_{q^{2}}^*$ with order $q^2-1$ and set $(u_0+v_0\a)^{q-1}=u+v\a$.
%  Then the correspondence matrix $\left(
%\begin{array}{cc}
%u& v\th \\
%v&u \\
%\end{array}
%\right)$ is an element of order $q+1$ of  $\SL(2,q)$ with order $q+1$, so  that  $\overline{{\small \left(
%\begin{array}{cc}
%u&\theta v \\
%v&u \\
%\end{array}
%\right)}}$ is an element of order $\frac{q+1}{2}$ in $\PSL(2,q)$.

Up to conjugacy,  set
$$M=\lg \overline{{\small \left(
\begin{array}{cc}
u&\theta v \\
v&u \\
\end{array}
\right)}},  \overline{ {\small \left(
\begin{array}{cc}
u'&-v'\theta \\
v'&-u' \\
\end{array}
\right)}}\di u, v, u', v'\in \FF_q, u^2-v^2\theta=1, -u'^2+v'^2\theta=1\rg\cong \D_{q+1},$$
noting that all the elements of the first type generate the cyclic subgroup of order $\frac{q+1}2$ of $M$ and all the elements  of the second type are all the noncentral involutions of $M$. Moreover,   $M$ has a central involution  if and only if   $q\equiv 3 (\mod 4)$, which is $\overline{{\small \left(
\begin{array}{cc}
0&\theta y \\
y&0 \\
\end{array}
\right)}}$ where $y^2=-\th^{-1}$.
 By Proposition~\ref{2,3},   every involution $\olg$ in $G$ has the form in the lemma, while $\olg \not\in M$ if and only if either  $a, b+c\th\ne 0$ or $a=0$ but $b\pm c\th\ne 0$. Now assume that $\overline{g}\in G\setminus M$.
 We shall  find the number of  noncentral involutions  $\olm\in M$ such that $|\olm\olg|>2$ and  $\lg\olm, \olg\rg$ is contained in a group which is isomorphic to $\D_{q+1}$, equivalently,
    such that $|\olm\olg|\di \frac {q+1}2$ and $|\olm\olg |>2$.
Now every  noncentral involution  $\olm\in M$ can be written as $m={\small \left(
\begin{array}{cc}
u_1&-v_1\theta \\
v_1&-u_1 \\
\end{array}
\right)}$,
where
\begin{eqnarray}
\label{uv}
-u_1^2+v_1^2\theta=1.
\end{eqnarray}  Then we have
$$mg={\small \left(\begin{array}{cc}
u_1a-v_1c\theta&u_1b+v_1a\theta\\
v_1a-u_1c&v_1b+u_1a\\
\end{array}
\right)}.$$
Now, $\det(\ld E_2-mg)=\ld^2-x\ld +1$ where
\begin{eqnarray}
\label{x}
x=\Tr (mg)=2au_1+v_1(b-c\theta).
\end{eqnarray}
Now $|\overline{mg}|\di \frac{q+1}2$ if and only if $\olm \olg$ has no fixed points
in $\PG(1,q)$, that is
\begin{eqnarray}
\label{x^2-1}
x^2-4\in N.
\end{eqnarray}
Note that $|\olm\olg|>2$ if and only if $x=\Tr(mg)\ne 0$. Combining Eq~(\ref{uv}) and Eq~(\ref{x}), we get
\begin{eqnarray}
\label{x-2}
((c\theta-b)^2-4a^2\theta)v_1^2+2x(c\theta-b)v_1+(x^2+4a^2)=0.
\end{eqnarray}
For Eq~(\ref{x-2}), since  $$4x^2(c\theta-b)^2-4[(c\theta-b)^2-4a^2\theta](x^2+4a^2)=-16a^2(-\theta x^2-t')=16a^2(\theta x^2+t')$$
where $t'=4\theta a^2-(c\theta-b)^2\ne 0$, it follows that
Eq(\ref{x-2}) has  two distinct  solutions for $v_1$ if and only if
\begin{eqnarray}
\label{x-4}
  x^2\in \frac{-t'}{\theta}+N.
\end{eqnarray}
Therefore,   $|\olm\olg|\di \frac{q+1}2$ and $|\olm\olg|>2$ if and only if
$x^2\in (4+N)\cap (\frac{-t'}{\theta}+N)\cap S^*$, that is, $(\frac x2)^2\in (1+N)\cap (\frac{-t'}{4\theta}+N)\cap S^*$.
Since $t'\ne 0$, it follows that $\frac{-t'}{4\theta}\ne 0$. Noting that  $\frac{-t'}{4\theta}=-a^2+\frac {1}{4\theta}(c\theta-b)^2=1+\frac {1}{4\theta}(c\th+b)^2$,
 we know that  $\frac{-t'}{4\theta}=1$  if and only if $b+c\theta=0$ that is, $\overline{g} \in M$, a contradiction.
 Therefore, $\frac{-t'}{4\theta}\neq 0, 1$.  Then by Lemma~\ref{feng},  we have at least $\lceil\frac{1}{8}(q-2\sqrt{q}-7)\rceil$  such   solutions  for $(\frac x2)^2$.
Furthermore, for each given $x$, we   have two solutions  for  $\olm$. Therefore, we get  at least  $2\lceil\frac{1}{8}(q-2\sqrt{q}-7)\rceil$    such noncentral involutions  $\olm $.
\qed

 \vskip 3mm

% (ii) $x=0$: Note that $x=0$ if and only  if  $t\in S$. If  $q\equiv 1(\mod 4)$, then  $-4\in N$ and so there is no solution for $x$.
%So $q\equiv 3(\mod 4)$.  In this case, $v=\pm \frac{2a}{\sqrt{t}}$ and $u=\mp \frac {b+c}{\sqrt{t}}$. Therefore,
 %  we get   only one the desired  element  $\olh$, while
% $\lg \olh, \olg \rg =D_4$.\qed

%\vskip 3mm

 Let $I=\{ a, b \cdots \}$ and $D=\{ \a, \b \cdots\} $ be the set of involutions    and  subgroups  isomorphic to $D_{q+1}$ of $G$, respectively.
 Note that for each $\a\in D$, there exists the unique  $z_\a\in I$ such that $\a=\C_G(z_\a)$, the centralizer  of $z_\a$ in $G$ for $q\equiv 3 (\mod 4)$.

To analyse the Saxl graph of $G$, we introduce a new graph $Y$, which is  the  bipartite  graph with partition $V=I\cup D$ where $\{ a, \a \}$ is an edge if $a\in \a$.
 Fix $\a\in D$ and let $Y_i(\a)$ denote the set of vertices $v$ in $V$  such that $d(\a, v)=i$, where $d(\a, v)=i$ means the distance between $\a$ and $v$ in graph $Y$ is $i$.
 First, we have the following lemma.
\begin{lem} \label{Gamma} The graph $Y$ is connected and of  diameter $4$ for $q\geq 17$ where $q$ is odd.
\end{lem}
\demo As above, let $\a:=M\cong \D_{q+1}$. Now $Y_1(\a)$ is the set of all the involutions of $M$.  By Lemma~\ref{connected}, for  any involution
$\overline{g}\in G\setminus M$, there exists an involution $\overline{m}\in M$ such that $|\overline{mg}|\di \frac{q+1}2$. Therefore, $\{\overline{m}, \overline{g}\}$ are contained in some subgroup isomorphic to  $\D_{q+1}$, which is a vertex in $Y_2(\a)$. In other words, $\overline{g}$ is adjacent to some vertex  in $Y_2(\a)$. This forces that for any $b\in I\setminus Y_1(\a)$, we have $b\in Y_3(\a)$, that is $d(Y)\le 4$.
Clearly, $d(Y)\ne 2$ and so $d(Y)=4$.\qed

\vskip 3mm

Now we address some results about some elements and subgroups in $\PSL(2,q)$ where $q$ is odd without proof. And then determine the lower bound of $|Y_2(\b)\cap Y_2(\a)|$.

\begin{lem} For $\PSL(2,q)$ where $q$ is odd,  the following statements  hold:
\begin{enumerate}
\item[\rm(i)] $\a\cap\b=1$ or $ \ZZ_2 $ (resp. $1,\ZZ_2$ or $\D_4$) for $q\equiv 1(\mod 4)$ (resp. $q\equiv 3(\mod 4)$).
\item[\rm(ii)] Every  subgroup $\ZZ_2$ is contained in  $\frac{q-1}{2} $ (resp. $ \frac{q+3}{2}$) subgroups $\D_{q+1}$ for $q\equiv 1(\mod 4)$ (resp. $q\equiv 3(\mod 4)$).
\item[\rm(iii)] There are $\frac{q(q+1)}{2}$ (resp. $\frac{q(q-1)}{2}$) subgroups $\ZZ_2$ in $\PSL(2,q)$ for $q\equiv 1(\mod 4)$ (resp. $q\equiv 3(\mod 4)$).
\end{enumerate} \end{lem}

\begin{lem} \label{length}   For $q\ge 17$, set  $w=w(q)=2\lceil\frac{1}{8}(q-2\sqrt{q}-7)\rceil$  and $n(\b)=|Y_2(\b)\cap Y_2(\a)|$. Then
\begin{enumerate}
\item[\rm(i)] $q\equiv 1(\mod 4)$, $q\ge 17$:   $n(\b)\ge \frac{q-1}2w-2$   or $\frac{q+1}2w$ for $\b\in Y_2(\a)$   and   $\b\in Y_4(\a)$, respectively;
\item[\rm(ii)] $q\equiv 3(\mod 4)$, $q\ge 17$: for any $\b$ such that $\b\cap \a=1$, $\b\cap \a=Z_2$ and $\b\cap \a=\D_4$, respectively, we have
$$n(\b)\ge  \frac{q+1}4w,  \frac{q-3}4w+\frac{q+1}4\,~  {\rm  and}\,~
\frac{q-3}4w+\frac{q+1}{4}.$$
\end{enumerate}
In particular,   $|Y_2(\b)\cap Y_2(\a)|\ge\frac{q-3}4w$ for all cases.

\end{lem}

\demo
{\it Case  $q\equiv 1(\mod 4)$.}

\vskip 3mm Now  the valency of $\a$ in $Y$ is $\frac{q+1}2$ and   for any $\b$, we have $\b\cap \a$ is either 1 or  $Z_2$. So there exist no cycles of length 4 in $Y$.
  Then
$|Y_2(\a)|=\frac{q+1}2\frac{q-3}2=\frac{q^2-2q-3}4$ and so
$|Y_4(\a)|= |\O|-1-\frac{q^2-2q-3}4=\frac{q^2-1}4.$

 Let  $\b\in Y_2(\a)$ and  set $\a\cap \b=\lg a_1\rg $.
  Note two facts: (i)
    there are $\frac{q-1}{2}-2$ vertices  $\b'$ in $Y_2(\a)\setminus \{\b \}$ adjacent to $a_1$;  (ii) by Lemma~\ref{connected}, for any involution $\overline{g}\in Y_1(\b)\setminus\{a_1\}$,  there are at least $w$ involutions $\overline{m}$ in $M$ such that  $|\overline{mg}|\di \frac{q+1}2$ and $|\overline{mg}|>2$, including $a_1$.
   Therefore,
   $$n(\b)\geq (w-1)(\frac{q+1}{2}-1)+(\frac{q-1}{2}-2)=\frac{q-1}{2}w-2.$$

Let   $\b\in Y_4(\a)$. Clearly,   $|Y_2(\a)\cap Y_2(\b)|\geq \frac{q+1}2w$.
\vskip3mm
{\it Case $q\equiv 3(\mod 4)$.}
 \vskip 3mm
 We prove it by the  following five steps.
 \vskip 3mm
  (1)  Note that  $\a\cong \D_{q+1}$. If   $|\a\cap \b|>4$, then there exists an element $g\in \a\cap \b$ where $|g|\di \frac{q+1}2$ and $|g|>2$. Then $\a=\N_G(g)=\b$.
So  for any distinct $\a, \b\in D$, we have $\a \cap \b$ is $\D_4$, $\ZZ_2$ or $1$. In  other words, in the graph  $Y$, we have   $|Y_1(\a)\cap Y_1(\b)|=0,1$ or 3 for any $\a, \b$.
For any $\b\in D$, by $z_\b$ we denote its center.

 If  $z_\b\in \a\cap \b$, then   $z_\a\in \a\cap \b$ and so
    $\{z_\a, z_\b, z_\a z_\b\}\subset \a\cap \b$, that is  $\a\cap \b=\D_4$.

Let $U=\{\C_G(a)\di a\in I\cap (\a\setminus \{z_\a\})\}$ which is contained in $Y_2(\a)$. Set $U'=Y_2(\a)\setminus U$.
Then for any $\b\in U$, there exists $b\in \a\setminus \{z_\a\}$ such that $\b=\C_G(b)$ and so  $|Y_1(\b)\cap Y_1(\a)|=3$. And for any $\b\in U'$, $|Y_1(\b)\cap Y_1(\a)|=1.$

 Moreover, for any $\b\in U$ where $\b=\C_G(b)$ for $b\in \a\setminus \{z_\a\}$,   let  $B=\{ b, z_\a, bz_\a\}\cup \{  \C_G(b), \C_G(z_\a), \C_G(bz_\a)\}$. Then
 the induced subgraph $Y(B)\cong K_{3,3}$.

  \vskip 3mm

(2)  Note $Y_1(\a)=I\cap \a$ with size $\frac{q+3}2$. Considering the number of edges between $Y_1(\a)$ and $Y_2(\a)$, we have
 $$\frac{q+3}2\cdot \frac{q+1}2=3|U|+|U'|=\frac{3(q+1)}2+|U'|,$$
which implies
$$|U'|=\frac{q^2-2q-3}4.$$
Therefore,   $$|Y_2(\a)|=|U|+|U'|=\frac{q+1}2+\frac{q+1}2\cdot \frac{q-3}2=\frac{q^2-1}4.$$
Since the diameter of $Y$  is 4,  except for the involutions in   $Y_1(\a)$,  all other involutions form the set $Y_3(\a)$ and so
$$ |Y_3(\a)|=\frac{q(q-1)}2-\frac{q+3}2=\frac{q^2-2q-3}2,$$ which forces
 $$|Y_4(\a)|=\frac{q(q-1)}2-1-\frac{q^2-1}4=\frac{q^2-2q-3}4.$$
 Since $G$ acts transitively on $I$ and $D$ by conjugacy,
   $\a $ can be any vertex  in $D$.

\vskip 3mm

(3) Case $\b\cap \a=1$:    Now $\b\in Y_4(\a)$. Consider a set
 $$W=\{ (a,b)\di a\in \a\setminus\{z_{\a}\}, b\in \b\setminus\{z_{\b}\},  |ab|>2, \, {\rm where}\, ~ a, b\in \g, \, {\rm for\, ~some }\, \g\in Y_2(\a) \}.$$
By Lemma~\ref{connected}, for any involution $b\in Y_3(\a)$,  there are at least
 $w$ nonentral involutions $a$ in $Y_1(\a)$ such that $|ab|\di \frac{q+1}2$ and  $|ab|>2$.  Therefore $|W|\ge \frac{q+1}2w$.

 Remind that in graph $Y$, for any 4-cycle: $\d_1, d_1, \d_2, d_2, \d_1$, at least one of $d_1$ and $d_2$ is a center of either $\d_1$ or $\d_2$. For any $\g\in Y_2(\a)\cap Y_2(\b')$, where $\b'\in Y_2(\a)$ or $Y_4(\a)$,  we have four possibilities:
$$\g\cap \a=\D_4,\,    \g\cap \b'=Z_2;\quad  \g\cap \a=Z_2, \, \g\cap \b'=\D_4;$$ $$\quad   \g\cap \a=Z_2, \,  \g\cap \b'=Z_2;\quad \g\cap \a=\D_4, \,  \g\cap \b'=\D_4,$$
 depending on  $z_\g\in Y_1(\a)\setminus Y_1(\b')$,  $z_\g\in Y_1(\b')\setminus Y_1(\a)$, $z_\g \not\in Y_1(\a)\cup Y_1(\b')$ and $z_\g \in Y_1(\a)\cap Y_1(\b')$, respectively. Therefore, from the value of $|W|$, we have  $n(\b)=|Y_2(\b)\cap Y_2(\a)|\ge \frac{q+1}4w$.

\vskip 3mm
(4) Case $\b\cap \a=Z_2$:
Write  $\b\cap \a=\{1, z_{\b_1}\}$ and $\b_2=\C(z_{\b_2})$ where $z_{\b_2}=z_{\b_1} z_{\b}$. Set
$$W=\{ (a,b)\di a\in \a\setminus\{z_{\a}\}, b\in \b\setminus \{z_{\b_1}, z_{\b_2}, z_{\b}\},  |ab|>2,  \, {\rm where}\,~ a, b\in \g, \, {\rm for\,~ some }\, \g\in Y_2(\a) \}.$$
As in (3), by Lemma~\ref{connected},  $|W|\ge \frac{q-3}2w$.
For any  $b\in \b\setminus \{z_{\b_1}, z_{\b_2}, z_{\b}\}$, if $b, z_{\b_1}\in \g$, where $\g\ne \b$,
then $\g=\b_1$ or $\g=\b_2$, which implies $|bz_{\b_1}|=2$, a contradiction. So $Y_1(z_{\b_1})\cap \{ Y_1(b)\mid b\in\{Y_1(\b)\setminus \{z_{\b_1}, z_{\b_2}, z_{\b}\} \}=\{\b\}$.
Therefore, as in (3),  from the value of $|W|$, we have
$$n(\b)=|Y_2(\b)\cap Y_2(\a)|\ge (\frac{q+3}2-2)+ \frac 12\cdot \frac{q-3}2(w-1)=\frac {q-3}4w+\frac{q+1}4.$$

\vskip 3mm (5) Case $\b\cap \a=\D_4$:   Then $\b\cap \a=\{1, z_\a, z_\b, z_\a z_\b\}$. Set $\g=\C_G(z_{\a}z_{\b})$ and $z_{\g}=z_\a z_\b$.
 Noting  a fact that any two vertices in $U$ intersect only at $Y_1(\a)$,   for any $b\in Y_1(\b)\cap Y_3(\a)$, we have $(Y_1(b)\cap Y_2(\a))\setminus \{\b\}\subset U'$.
Set
$$W=\{ (a,b)\di a\in \a\setminus\{z_{\a}\}, b\in \b\setminus \{z_{\a}, z_{\b}, z_{\g}\},  |ab|>2,\,   {\rm where}\, ~a, b\in \d, \, {\rm for\,~ some }\, \d\in Y_2(\a) \}.$$
As in (3), by Lemma~\ref{connected}, we have $|W|\ge \frac{q-3}2w$. Therefore, $|\{Y_1(b)\cap Y_2(\a)\di b\in Y_1(\b)\cap Y_3(\a)\}|\ge \frac{q-3}{4}(w-1).$
By adding  the vertices $U\setminus \{\b\}$,  we get
 $$n(\b)=|Y_2(\a)\cap Y_2(\b)|\ge \frac{q-3}4(w-1)+\frac{q-1}{2}=\frac{q-3}{4}w+\frac{q+1}{4}.$$
\qed

\subsubsection{ $G=\PSigmaL(2,q)$, $M=\D_{q+1}.\ZZ_n$} \label{PSigmaLPGammaL}
Let  $p$ be an odd prime and $q=p^n\ne 7,9.$
Let $G=\PSigmaL(2,p^{n})$ and $M=\D_{q+1}.\lg fu\rg $ for some  $u\in \PSL(2,q)$, where $f$ is induced by  the field automorphism  $\sigma: a\mapsto a^p$ for any $a\in \FF_{p^{n}}^*$.
Let $T:=\PSL(2,p^{n})$ and $M_0=\D_{q+1}$.  Set $\a=M$, $\O'=[T:M_0]$ and $\O=[G:M]$.  Consider   the primitive permutation representation of $G$ on $\O$ with right coset multiplication action. Then the action of $T$ on $\O$ is
equivalent to that on $\O'$.

 Let $\FF_{p^{n}}^*=\lg \d\rg $ and  $\FF_{p^{2n}}=\FF_{p^n}(\zeta)$, where $\zeta^2=\th :=\d^{\frac{p^n-1}{2^l}}$ and $p^n-1=2^lk$  for an odd integer $k$. In particular, $|\th |=2^l$.   Take a generator $x+y\zeta$ of $\FF_{p^{2n}}^*$. Set
$$s={{ \overline{\left(
\begin{array}{cc}
x&y\th \\
y&x \\
\end{array}
\right)  }}},\quad  t={{ \overline{\left(
\begin{array}{cc}
1&0 \\
0&-1 \\
\end{array}
\right)  }}} \quad {\rm and}\quad w={{ \overline{\left(
\begin{array}{cc}
1&0 \\
0&\th^{\frac{1-p}2} \\
\end{array}
\right)  }}} .$$

\begin{lem}   \label{M} Up to conjugacy in $G$,  we have
 $$ M_0=\lg a, b\rg =\D_{q+1},\quad  and \quad  M=N_G(M_0)=M_0\lg c\rg =\D_{q+1}.Z_n,$$
  where $a=s^2$ and

  \begin{enumerate}
\item[{\rm (i)}]  $p\equiv 3(\mod 4),$ $n$ is odd:  $b=st, c=fw$;
\item[{\rm (ii)}]   $p\equiv 1(\mod 4):$  $b=t, c=fw;$
\item[{\rm (iii)}]  $p\equiv 3(\mod 4),$  $2\di n: b=t, c=fws^{\frac{q+1}2}.$\end{enumerate}
 Moreover,  we have
 \begin{enumerate}
\item[{\rm (iv)}] $a^c=a^{p}$;
\item[{\rm (v)}]  $c^{\frac{n}{r}}=f^{\frac{n}{r}}t$, for any  odd prime divisor $r$ of $n$; $c^{\frac{2n}{r}}=f^{\frac{2n}{r}}$;  $c^n=t$.
\end{enumerate}
%For cases (ii) and (iii) we have $t^c=t$.
\end{lem}
\demo Since $|s|=q+1$ and $T$ has only one conjugacy class of  subgroups of order $q+1$, we take $a=s^2$. Note that
   in what follows, we find  $b\in T$  and $c\in G$ such that  $\lg a, b\rg =\D_{q+1}$ and $c$ normalizes $\lg a, b\rg $.

If $p\equiv 3(\mod 4)$ and $n$ is odd, then  $t\in \PGL(2,q)\setminus T$ and so let  $b=st$ so that $b\in T$. Moreover,   take $c=fw$ so that $c$  normalizes $\lg a, b\rg $.

If $p\equiv 1(\mod 4)$,  then $t, w\in T$  and we set $b=t$ and  $c=fw$. One can check that $c$ normalizes $\lg a, b\rg $.

 If $p\equiv 3(\mod 4)$  and  $2\di n$, then $t\in T$  but $w\in \PGL(2,q)\setminus T$. So we take $b=t$ and $c=fws^{\frac{q+1}2}$.

Note that
there is an injective homomorphism $\phi$: $\FF_{q^2}^*/\FF_q^*\to \PGL(2,q)$, given by $(x+y\zeta)\FF_q^* \longmapsto {{\overline{\left(
\begin{array}{cc}
x&y\theta \\
y&x\\
\end{array}
\right)  }}}$. So
$$s^p=\phi((x+y\zeta)^p\FF_q^*)=\phi((x^p+y^p\theta^{\frac{p-1}{2}}\zeta)\FF_q^*)=
{{ \overline{\left(
\begin{array}{cc}
x^p&y^p\theta^{\frac{p-1}{2}}\theta \\
y^p\theta^{\frac{p-1}{2}}&x^p \\
\end{array}
\right)  }}}$$ for all the three cases. Moreover,  $$s^f={{ \overline{\left(
\begin{array}{cc}
x^p&y^p\theta^{p} \\
y^p&x^p \\
\end{array}
\right)  }}} \,\quad {\rm and}\, \quad s^{fw}={{ \overline{\left(
\begin{array}{cc}
x^p&\theta^{\frac{p-1}{2}}y^p\theta \\
\theta^{\frac{p-1}{2}}y^p&x^p \\
\end{array}
\right)  }}}.$$
 Then  $s^c=s^p$.

  Suppose that either $p\equiv 3(\mod 4)$ for $n$ is odd or  $p\equiv 1(\mod 4)$ and $r$ is an odd prime divisor of $n$. Then $$c^{\frac{n}{r}}=(fw)^{\frac{n}{r}}=f^{\frac{n}{r}}w^{f^{\frac{n}{r}-1}}w^{f^{\frac{n}{r}-2}}\cdots w^fw.$$ That is to say $$c^{\frac{n}{r}}=f^{\frac{n}{r}}{{ \overline{\left(
\begin{array}{cc}
1&0 \\
0&(\th^{\frac{1-p}2})^{\frac{1-p^{\frac{n}{r}}}{1-p}} \\
\end{array}
\right)  }}} =f^{\frac{n}{r}}{{ \overline{\left(
\begin{array}{cc}
1&0 \\
0&-1 \\
\end{array}
\right)  }}}=f^{\frac{n}{r}}t\quad{\rm and}\quad
c^n=f^n{{ \overline{\left(
\begin{array}{cc}
1&0 \\
0&(\th^{\frac{1-p}2})^{\frac{1-p^n}{1-p}} \\
\end{array}
\right)  }}}=t.$$

If $p\equiv 3(\mod 4)$ for $n$ is even and  $r$ is an odd prime divisor of $n$.  Set $$h=ws^{\frac{q+1}{2}}={{ \overline{\left(
\begin{array}{cc}
0&\theta z\\
\theta^{\frac{1-p}{2}}z&0 \\
\end{array}
\right)  }}} \quad {\rm and}\quad s^{\frac{q+1}{2}}={{ \overline{\left(
\begin{array}{cc}
0&z\theta \\
z&0 \\
\end{array}
\right)  }}} $$ for some $z\in \FF_q^*$.
Then $$c^2=(fh)^2=f^2h^fh=f^2{{ \overline{\left(
\begin{array}{cc}
(\theta^{\frac{1-p}{2}})\theta^pz^{1+p}&0 \\
0&(\theta^{\frac{1-p}{2}})^p\theta z^{1+p} \\
\end{array}
\right)  }}} =f^2{{ \overline{\left(
\begin{array}{cc}
1&0 \\
0&\theta^{\frac{1-p^2}{2}} \\
\end{array}
\right)  }}}.$$
It follows that
$$c^{\frac{n}{r}}=(c^2)^{\frac{n}{2r}}=f^{\frac{n}{r}}{{ \overline{\left(
\begin{array}{cc}
1&0 \\
0&\theta^{\frac{1-p^{\frac{n}{r}}}{2}} \\
\end{array}
\right)  }}}=f^\frac{n}{r}t \quad {\rm and}\quad c^n=(c^2)^{\frac{n}{2}}=f^n{{ \overline{\left(
\begin{array}{cc}
1&0 \\
0&\theta^{\frac{1-p^n}{2}} \\
\end{array}
\right)  }}}=t.$$
\qed

\vskip 3mm
We  first need to know which regular suborbits $\a^{gM_0}$  of $T$  is fixed by a nontrivial subgroup of $\lg cu\rg$.
Suppose that $\a^{gM_0}$  is such a regular suborbit.  Then it must be fixed by an element $c'u_0$ of prime order $r'$ for some $u_0\in M_0$, where
$c'=c^{\frac {|c|}{r'}}$.  Then
the point $\a^{g}$ is fixed by a subgroup $K$ of order $r'$, where $K\leq M$ and $K\cap M_0=1$.
In what follows, let $n'(r')$ be the number of regular suborbits of $T$ which are fixed by $K$.
The following lemma determines the upper bound of $n'(r')$.
\begin{lem} \label{conjugacy} With the above notations, we get the following results.
\begin{enumerate}
\item[{\rm (1)}]   If   a  subgroup $K$ in  $M\setminus M_0$  has order a prime $r'$, then $r'$ is odd.
\item[{\rm (2)}]  All the subgroups of order $r'$ in $M\setminus M_0$ are conjugate to $\lg f'\rg $ in $M$, where $f'=f^{\frac{n}{r'}}$.
\item[{\rm (3)}]  $n'(r')\le [\frac{1}{2}p^{\frac{n}{r'}}(p^{\frac{n}{r'}}-1)/(p^{\frac{n}{r'}}+1)]$.
\end{enumerate}
\end{lem}

\demo (1)  Suppose that $K\leq M$ and $K\cap M_0=1$ with  prime order $r'$. Suppose that $r'=2$.  Then $n$ is even and we are in the cases (ii) and (iii) in Lemma~\ref{M} and in particular $c^n=t$.
  Then $K=\lg uc^l\rg \cong Z_2,$ where $u=a^it^j\in M_0$.
Then $1=(a^it^jc^l)^2=c^{2l}a^{ip^{2l}+i(-1)^jp^l}$, which implies $c^{2l}=1$, that is $c^l\in \lg t\rg$ and $K\le M_0$, a contradiction.
  Therefore,  $r'$ is odd.

\vskip 3mm
(2)  Let $n=r'^ln'$, where $r'$ is odd and $(r', n')=1$.  Let $K$ be a subgroup of order $r'$ in $M$. If $r'\nmid |M_0|$, then  we take $K=\lg f'\rg $ up to conjugacy where $f'=f^{\frac{n}{r'}}$.

Suppose that $r'\di |M_0|$, that is $r'\di \frac{p^n+1}2$. Set $\frac{p^n+1}2=r'^ji$, where $(r', i)=1$.

Consider a Sylow $r'$-subgroup $\lg a^i\rg \lg c^{\frac{2n}{r'}}\rg$  in $M_0\lg c^{\frac{2n}{r'}}\rg $.
If  $$(a^i)^{c^{\frac{2n}{r'}}}a^{-i}=a^{ip^{\frac{2n}{r'}} } a^{-i}=a^{i(p^{\frac{2n}{r'}}-1)}=1,$$
 then $r'^j\mid p^{\frac{2n}{r'}}-1$ and so  $r'^j\mid p^{\frac{n}{r'}}+1$.
 Set  $p^{\frac{n}{r'}}=kr'^j-1$
 for some integer k. Then
  $$p^n+1=(kr'^j-1)^{r'}+1\equiv 0(\mod r'^{j+1}),$$
  a contradiction.  Therefore,
  $\lg a^i\rg : \lg c^{\frac{2n}{r'}}\rg $  is an nonabelian Sylow $r'$-subgroup, which   contains a maximal cyclic subgroup.
  By Lemma~\ref{maxcyclic},  we may set $K=\lg c^{\frac {2n}{r'}}\rg =\lg f^{\frac {2n}{r'}}\rg $ so that   we have $K=\lg f'\rg $, where $f'=f^{\frac n{r'}}$.
\vskip 3mm
(3) Let $K=\lg f'\rg $.  Since $$|\Fix (K) |=\frac {|\N_G(K)|}{|\N_M(K)|}= \frac {|\C_G(K)|}{|\C_M(K)|}=\frac{p^{\frac{n}{r'}}(p^{\frac{n}{r'}}-1)}{2}$$
 and $|\N_{M_0K}(K):K|= p^{\frac{n}{r'}}+1$, we have $f'$ fixes   $$|n'(r')|\le \frac{1}{2}p^{\frac{n}{r'}}(p^{\frac{n}{r'}}-1)/(p^{\frac{n}{r'}}+1)$$
regular suborbits of $T$.
\qed

\vskip 3mm
\begin{lem} \label{PSigmaL} The diameter of the Saxl graph $\si(G)$ is $2$.\end{lem}
\demo
If $q\in\{5,11,13\}$, then checked by Magma, we have that the Saxl graphs have diameter $2$.

Next we suppose that $q\geq 17$.

Let $n'$ be the number of regular suborbits of $T$ which are fixed setwise  by a
nontrivial element of $\lg c\rg$. In particular, $n'=0$ if $f=1$, that is $G=T$.

Suppose that   $n=2^{e_0}n_0$, where $n_0=r_1^{e_1}r_2^{e_2}\cdots r_l^{e_l}$ with  $3\le r_1<r_2 \cdots <r_l$ are distinct primes.
 Also, $n_0\ge 5^{l-1}r_1$,
that is $l\le 1+\log_5\frac {n_0}{r_1}.$ Then we have
$$n'\le \sum_{i=1}^ln'(r_i)\le l\cdot{\rm max}\{n'(r_i)\di 1\le i\le l\}=ln'(r_1)
\le \frac 12(1+\log_5\frac {n_0}{r_1})p^{\frac {n}{r_1}}.$$

 Let $Z$ be the complement of $\si(G)$.  Considering the relation between respective suborbits relative to $M$ and $M_0$,  for
any $\b \in Z_1(\a)$, we have
 $$|Z_1(\a)|=|Y_2(\a)|+n'(q+1),\quad  \quad |Z_2(\a)|=|Y_4(\a)|-n'(q+1),$$
$$|Z_1(\b)\cap Z_2(\a)|\le |Y_2(\b)\cap Y_4(\a)|+n'(q+1)\le  |Y_2(\b)\setminus \{\a\}|-|Y_2(\b)\cap Y_2(\a)|+n'(q+1).$$
Note that $Z_1(\a)\cup Z_2(\a)\cup \{\a\}=V(\Sigma(G))$ and $Z_2(\a)=\Sigma_1(\a)$ where $\Sigma_1(\a)$ denotes the set of neighbors of $\a$ in graph $\Sigma(G)$.
Therefore, $\si(G)$ has diameter $2$ if
$|Z_1(\b)\cap Z_2(\a)|<|Z_2(\a)|$, but this  holds  provided
 $$|Y_2(\b)\setminus \{\a\}|-|Y_2(\b)\cap Y_2(\a)|+n'(q+1)<|Y_4(\a)|-n'(q+1),$$
 that is,
  $$\begin{array}{lll}&&2n'(q+1)<k:=|Y_4(\a)|-|Y_2(\b)\setminus \{\a\}|+|Y_2(\b)\cap Y_2(\a)|\\
  &=&\mp\frac{q^2-2q-3}4\pm\frac{q^2-1}4+\e+|Y_2(\b)\cap Y_2(\a)| \end{array}$$
 where $\e= 1$ (resp. 0) if $\b\in Y_2(\a)$ (resp. $Y_4(\a))$.
 Further, we have $k\ge \pm\frac{q+1}2+n(\b)$ for $ q\equiv \pm 1(\mod 4)$ where
 $n(\b)=|Y_2(\b)\cap Y_2(\a)|$.
 % Note that $k:=|Y_4(\a)|-(|Y_2(\b)\setminus \{\a\}|+|Y_2(\b)\cap Y_2(\a)|)=1+|Y_4(\a)|-(|Y_2(\b)|-|Y_2(\b)\cap Y_2(\a)|)$ for $\b\in Y_2(\a)$, we have $\b\geq 1$.
 Then we have the following possibilities:

\vskip 3mm
(i) $q=p\ge 17:$ In this case, $n'=0$ and we only need to consider $\b\in Y_2(\a)$.  By  Lemma~\ref{length}, $n(\b)\ge \frac{q-1}2w-2$ for $q \equiv 1 (\mod 4)$ and $n(\b)\ge \frac{q-3}4w+\frac{q+1}{4}$ for $q \equiv 3 (\mod 4)$. Then we get  $k>0$.

\vskip 3mm
(ii) $q=p^n\ge 25$ and $n\ge 2$: By Lemma~\ref{length} again, we get that
$$k\ge \frac{q-3}4w-\frac{q+1}2\ge \frac 1{16}(q^2-2q^{\frac 32}-18q+6q^{\frac 12}+13)\ge \frac 1{16}q^2(1-2q^{-\frac 12}-18q^{-1}).$$
By Lemma~\ref{conjugacy}.(1) and (3), for $n=2^e$,  $n'=0$;  for  an odd prime $n$, $n'=0$, $n'\le 1$ and $n'\le 2$, if $p=3, 5, 7$, respectively.
  One may check that  for these cases, we have $2n'(q+1)<k$. So in what follows, we let either  $n=3$ and $p\ge 11$ or $n\ge 5$.
Since
$$(1+\log_5\frac {n_0}{r_1})p^{\frac {n}{r_1}}(q+1)< (1+\log_5n)q^{\frac 13}q(1+3^{-6})=3^{-6}\cdot 730\cdot (1+\log_5n)q^{\frac 43},  $$
$$\frac 1{16}q^2(1-2q^{-\frac 12}-18q^{-1})\ge
\frac 1{16}q^2(1-2\cdot 3^{-3}-18\cdot 3^{-6})\ge \frac 1{16}\times 657\times 3^{-6}q^2 ,$$
it suffices to show
$18(1+\log_5n)\le p^{\frac {2n}3}.$ If $n=3$ and $p\ge 11$, this inequality holds. So assume $n\ge 5$. Then it suffices to show $18(1+\log_5n)\le 3^{\frac {2n}3},$
which is clearly true for $n\ge 5$.

In summary,  the Saxl graph has diameter $2$ for any $q$ provided $b(G)=2$.
\qed

\subsubsection{ Proof of Theorem~\ref{main1} for $M_0=D_{\frac{2(q+1)}{d}}$ where $d=(2,q-1)$}
The  first two lemmas below show $b(G)>2$ for either even $q$ or $G=\PGL(2,q)$.

\begin{lem}\label{q=even}
 With our notations, suppose that $q=2^n\ge 8$.
Then $b(G)>2$ for $T\leq G\leq \PGammaL(2,q)$.
\end{lem}
\demo  Suppose that $q=2^n$. It suffices to  show $b(G)>2$ when $G=\PSL(2,q)$ and $M=\D_{2(q+1)}.$  In the group $G$, any two subgroups isomorphic to $\D_{2(q+1)}$ intersect either in $\ZZ_2$ or in $1$. Let $K$ be a subgroup of order 2 in  $G$. Then there is only one conjugacy class isomorphic to $\ZZ_2$ in $\D_{2(q+1)}$. Since $\N_G(K)\cong \ZZ_2^n$ and $\N_M(K)=K$, $K$ fixes $2^{n-1}$  points in $\O$.
In each suborbit of length $q+1$, $K$ fixes one point and so $G$ has $2^{n-1}-1$ such suborbits.
Observing  $(q+1)(\frac{q}{2}-1)+1=|\O|$, we know that there exists no regular suborbit, that is $b(G)>2.$  \qed

\begin{lem}\label{dq+1}
If  $G=\PGL(2,q)$, then $b(G)>2$.
\end{lem}
\demo By Lemma \ref{q=even}, the lemma is true for even $q$. Next, suppose that $q$ is odd. Let $M_1,M_2\in G$ be two distinct subgroups which are isomorphic to $\D_{2(q+1)}$. Then there exist two distinct involutions $a_1,a_2\in G$ such that $\C_G(a_i)=M_i$ and $\lg a_1,a_2 \rg\cong \D_{2k}$ where   either $k=p$ or $k$ is a divisor of $q\pm 1$. Note that $a_i$ has no fixed points in $\PG(1,q)$ for $i=1,2$.
Since every element in $\D_{2p}$ has fixed points,  it follows that $k\neq p$.  If $k \di q\pm1$, then either $\D_{2k}\leq \D_{2(q+1)}=\C_G(a_3)$ or $\D_{2k}\leq \D_{2(q-1)}=\C_G(a_4)$, for some involutions $a_3, a_4\in G$. Correspondingly,  either $a_3$ or  $a_4$ is included in $M_1\cap M_2$. This follows $b(G)>2$.
\qed

\vskip 3mm
Now we are ready to prove Theorem~\ref{main1}.
\vskip 3mm
\f {\bf Proof:}  By Lemma~\ref{q=even}, if $q$ is even, then $b(G)>2$ for any $T\le G\le \PGammaL(2,q)$.  So let $q$ be odd.
Using Lemma~\ref{PSigmaL}, we obtain that $b(G)=d(\Sigma(G))=2$ for $G=\PSigmaL(2,q)$.
By Lemma~\ref{subgraph},  $d(\si(G_1))=2$ for $G_1=T:\lg f^i\rg $ where $0\le i\le n-1$.

By Lemma~\ref{dq+1}, $b(G)>2$ if $G=\PGL(2,q)$.  We need to consider other  groups $T<G_2\le \PGammaL(2,q)$ but $ \PGL(2,q)\nleq G_2$.
If $G_2/T$ is not cyclic, then $G_2/T$ contains $\overline{\delta}$, which implies that $G_2$ contains  $\PGL(2,q)$ and so $b(G_2)\ne 2.$   Suppose that  $G_2/T$ is cyclic, that is $G_2/T=\lg \overline{\delta}\overline{f}^i\rg$. If $|f^i|$ is odd, then $\overline{\delta}\in G_2/T$ and $G_2$ contains $\PGL(2,q)$ again.
So $|\delta f^i|=|f^i|=2k$ for some integer $k$.  Let  $\Del$ be the union of all regular suborbits $\Del_j$ of  $T$  such that $\Del_j\cap \Del_j^{f^l}=\emptyset$ for any $1\leq l\leq n-1$. Then $\Del$ is the neighbourhood of $\a$ in Saxl graph $\si(\PSigmaL(2,q))$.
 For every such suborbit $\Del_j$,   we have $\Del_j^{\delta}=\Del_j$. Therefore,    $\Del_j\cap \Del_j^{(\delta f^i)^m}=\Del_j\cap \Del_j^{f^{im}}=\emptyset$ for any $m$ such that $n\nmid im$. Thus,   the neighbourhood of $\a$ in Saxl graph $\PSL(2,q).\lg \delta f^i\rg $ contains $\Del$.  So the diameter of Saxl graph for $\PSL(2,q).\lg \delta f^i\rg$ is also 2.\qed

\subsection{$M_0=\PSL(2, p^m)$}

The main result of this section is the following Theorem~\ref{main2}.
\begin{theorem}\label{main2}
Suppose that $T=\soc(G)=\PSL(2,q)$ and $M$  a maximal subgroup of $G$ with $M_0=M\cap T=\PSL(2,p^m)$,  where $p^m>2$,  either $\frac nm$  is  an odd prime or $p=2$ and $n=2m$. Consider the primitive permutation representation of $G$ on $\O=[G:M]$. Then  $b(G)=2$ if and only if $\frac nm$  is an odd prime. Moreover, the corresponding  Saxl graph $\si(G)$ has diameter $2$.

 \end{theorem}
\demo In Lemma~\ref{PSL}, we shall compute $|\G|$ and show $b(G)>2$ for $n=2m$ and $p=2$.   Suppose $\frac nm$ is odd.
Note that the action of $T$ on $\O=[G:M]$ is equivalent to that of $\O'=[T:M_0]$. We get that the regular suborbits of $T$ exist for an odd prime $\frac{n}{m}$ by Lemma~\ref{PSL}.
In Lemma~\ref{PGLPGL}, we shall show that  for any  $T$-regular suborbit $\Delta$, there exists another $T$-regular suborbit $\Delta'$  such that $\Delta \cup \Delta'$ is a $\PGL(2,q)$-regular suborbit.
 Finally, in Lemma~\ref{PGammaL}, we show that the Saxl graph of  $\PGammaL(2,q)$ has diameter 2. Therefore, by Lemma~\ref{subgraph}, for every group $G$ with socle $T$, its Saxl graph has diameter 2. This completes the proof.
 \qed

\subsubsection{$G=T, M=\PSL(2, p^m)$}
Let  $G=T=\PSL(2,q)$ where $q=p^n$ and $G_\a=M=\PSL(2,p^m)$ where $p^m>2$ and either $\frac nm$  is  an odd prime or $p=2$ and $n=2m$. Set $K=M_\b$ where $\b\in \O\setminus \{\a\}$ and let $\Fix(K)$ be
 the set of fixed points of $K$ in $\O$.
  Then we shall  determine the  suborbits of $G$ relative  to $\a$,  according to  the possibilities of the subgroups $K$ listed  in Proposition~\ref{psl}.

\begin{lem}\label{sl}
Let $G=\PSL(2,q)$ where $q=p^n$ and $G_\a=M=\PSL(2,p^m)$ where either $\frac nm$  is  an odd prime or $p=2$ and $n=2m$. Then the following statements hold.
\begin{enumerate}
\item[\rm(i)] Let $P$ be a  Sylow $p$-subgroup of $M$ and  $P_1$ a subgroup of order $p$ of $P$. Then  $\Fix(P)=\Fix(P_1)$.
\item[\rm(ii)] Let $L$ be a subgroup of $M$ with order  $\frac 1d(p^m\pm 1)$ where $d=(2,q-1)$ and $R$ a subgroup of prime order $r$ of $L$.  Then $\Fix(L)=\Fix(R)$.
\end{enumerate}
\end{lem}
\demo (1) Take a Sylow $p$-subgroup $P=\{ a(y):=\overline{{\small\left(
\begin{array}{cc}
1&0 \\
y&1 \\
\end{array}
\right)}} \di y\in \FF_{p^m}\}$ of $M$. Let $P_1=\lg a(x_1)\rg$ and $P_2=\lg a(x_2)\rg$ be any two subgroups in $P$ with order $p$ such that $P_1^g=P_2$ for some $g\in G$.
Suppose $g=\overline{{\small \left(
\begin{array}{cc}
a&b\\
c&d \\
\end{array}
\right)}} \in G$ such that $P_1^g=P_2$. Then
 there exist $y_1, y_2\in \FF_{p^m}^*$ such that $\overline{{\small \left(
\begin{array}{cc}
1&0 \\
y_1&1 \\
\end{array}
\right)}}^g=\overline{{\small \left(
\begin{array}{cc}
1&0 \\
y_2&1 \\
\end{array}
\right)}}$.
By solving this equation, we get $b=0$, $d=a^{-1}$ and $dy_2=ay_1$.
 Then $a^2\in \FF_{p^m}^*$.
Since $\frac{p^n-1}{p^m-1}$ is odd, we know  the quotient group $\FF_{p^n}^*/\FF_{p^m}^*$ is of odd order.  So $a^2\in \FF_{p^m}^*$ if and only if $a\in \FF_{p^m}^*.$
Taking $c=0$, we get that $g\in \PSL(2,p^m)$, that is $P_1$ and $P_2$ are conjugate in $\PSL(2,p^m)$.
Now $\N_M(P_i)=P:\lg z \rg$ and  $\N_G(P_i)=\ZZ_p^n:\lg z \rg$ for some $\lg z\rg \in M$. Therefore,  we get   $|\Fix(P_i)|=p^{n-m}$. Also, $\N_M(P)=P:\lg z'\rg$ and  $\N_G(P)=\ZZ_p^n:\lg z'\rg$ for some $\lg z'\rg \in M$, we have
$|\Fix(P)|=p^{n-m}$ too. This gives $\Fix(P_i)=\Fix(P)$.

\vskip 3mm
(2)   Suppose that $L\le M$ such that $L\cong \ZZ_{\frac 1d(p^m\pm 1)}$ and $R\le L$ with $|R|=r$, where $r$ is a prime divisor of $\frac 1d(p^m\pm 1)$.
Then  $\N_G(L)=\N_G(R)$ is isomorphic to $\D_{\frac{2(p^n\pm 1)}d}$ for $\frac{n}{m}$ is an odd prime and to  $\D_{2(p^{2m}- 1)}$  for $n=2m$;
 and $\N_M(L)=\N_M(R)\cong \D_{\frac{2(p^m\pm 1)}d}$. Note that there is only one conjugacy class of subgroups  of order  a divisor of $\frac{p^m\pm 1}d$ in   $M$. Therefore,
 $|\Fix(L)|=|\Fix(R)|=\frac{p^n\pm 1}{p^m\pm 1}$  or  $\frac{p^{2m}-1}{p^m\pm 1}$ for $\frac{n}{m}$ is an odd prime or $n=2m$, respectively. \qed

\vskip 3mm
{\it  From now on, in several lemmas, we shall use Proposition~\ref{man} to compute the size $|\Fix(K)|$ where $K$ is a subgroup of $M$. To do this (see the proposition), we first need to know the number of conjugacy classes in $M$ of all $G$-conjugate subgroups of $K$  which are contained in $M$.  From this information  about $|\Fix(K)|$, we may determine the number of suborbits whose point stabilizers in $M$ are $K$.}

\begin{lem} \label{D4} Suppose that $K$ is one of the following:
\begin{enumerate}
\item[\rm(i)] $\D_4$ where $p$ is odd;
\item[\rm(ii)] $\D_{2l}$ where $l$ is odd and $l\di \frac{p^m\pm 1}d$;
\item[\rm(iii)] $\ZZ_p^s:\ZZ_t$ where $s\leq m$,  $t>1$  and $t\mid \frac{p^s-1}{d}$.
\end{enumerate}
Then $\Fix(K)=\{ \a\}$.  So there exist no nontrivial suborbits $\b^M$ with $M_\b=K.$
\end{lem}
\demo Note that $G$   has one conjugacy class of  subgroups isomorphic to   $\ZZ_{\frac{p^n\pm 1}d}$ and $\D_{\frac{2(p^n\pm 1)}d}$, respectively,  where $d=(p^n-1,2)$.

 \vskip 3mm

 (1)  $K=\D_4$, where $p$ is odd.

 Let $N=\lg a, b\rg \cong \D_{\frac{2(p^n\pm 1)}d}$, where $|a|=\frac{p^n\pm 1}d$ and $|b|=2$. Set $a_1=a^e$, where $e=\frac{p^n\pm 1}{p^m\pm 1}$, which is  odd.  Set $z=1$ for either $p=2$ or  $p^n\pm 1\equiv 2(\mod 4)$ and $z=a^{\frac{(p^n\pm 1)}{4}}$ for $p^n\pm 1\equiv 0(\mod 4)$. Then $\Z(N)=\lg z\rg$.

\vskip 3mm

As stated before, in order to get the size of fixed points of $K$ by using Proposition~\ref{man},  we need to find the number of conjugacy classes in $M$ of $G$-conjugacy subgroups of $K$ which are contained in $M$.

Note that if $p^m\equiv 3, 5(\mod 8)$, then $M$ has only one conjugacy class of subgroups isomorphic to $\D_4$.

If $p^m\equiv 1, -1(\mod 8)$, then $M$ has  two conjuacy classes of subgroups isomorphic to $\D_4$ with respective  representatives $K_1=\lg z, b\rg$ and $K_2=\lg z, ba_1\rg$, noting that $K_1, K_2\le \lg a_1,b\rg \cong \D_{\frac{2(p^m\pm 1)}d}$. Suppose that there is an element $g\in G$ such that $K_1^g=K_2$. Since, for $i=1,2$  a subgroup isomorphic to  $A_4$ of $\N_G(K_i)$ is transitive on the three involutions in $K_i$, we may choose $g' \in\C_G(z)=N$ such that $\{b, bz\}^{g'}=\{ba_1, ba_1z\}$. For any  $g'=b^la^k$ where $l$ and $k$ are integers, we have $b^{g'}=ba^{2k}$, which is impossible.
 Therefore, $K_1$ and $K_2$  are not conjugate in $G$ either. This means that the two classes of subgroups isomorphic to $\D_4$ in $M$ cannot be  merged in $G$.

Since $A_4\leq \N_M(K)\leq \N_G(K)\leq S_4$ and $|G:M|$ is odd, $\N_G(K)=\N_M(K)$. By Proposition \ref{man}, it follows that $|\Fix(K)|=|\N_G(K):\N_M(K)|=1$, that is $\Fix(K)=\{\a\}$.

\vskip 3mm
(2)  $K=\D_{2l}$ where $l$ is odd and $l\di \frac{p^m\pm 1}d$.
   Then  up to conjugacy, set $K=\lg a_2,b\rg \le M$, where $a_2=a_1^{\frac{p^m \pm 1}{ld}}$.

Note that in $G$, the normalizer of a dihedral group $\D_{2l}$ is either  $\D_{2l}$ or $\D_{4l}.$
First, suppose  $p=2$. Then $l\di (2^m\pm 1)$. Since both $M$ and $G$ have one class of such groups and  $\N_G(K)=\N_M(K)=K$, we get   $|\Fix(K)|=1$, that is $\Fix(K)=\{\a\}$.
Next, suppose  that $p$ is odd.  Then  $M$ has only one conjugacy class of subgroups isomorphic to $\D_{2l}$.
Since $e$ is odd, $\frac{p^m\pm 1}{2}/l$ is even if and only if $\frac{p^n\pm 1}{2}/l$ is even.
Then $\N_G(K)=\N_M(K)$. Thus $|\Fix(K)|=1$.

 \vskip 3mm
(3)  $K=\ZZ_p^s:\ZZ_t$ where  $s\leq m$, $t>1$ and $t\mid \frac{p^s-1}{d}$.

Since $\ZZ_p^s$ is a characteristic subgroup of $K$, we have $\N_G(K)\leq \N_G(\ZZ_p^s)=\ZZ_p^n:\ZZ_{l}$ with a similar argument in the first part of Lemma \ref{sl}, where $l\di(p^m-1)$. For any $g\in\ZZ_p^n\cap \N_G(K)\leq \N_G(K)$, there is an element $h\in\ZZ_p^s$ such that $\ZZ_t^g=\ZZ_t^h$. Thus $gh^{-1}\in \ZZ_p^n\cap \N_G(\ZZ_t)=1,$ noting  $\ZZ_p^n:\ZZ_t$ is a  Frobenius group.  It follows that $g=h\in \ZZ_p^s$.  Thus we have $\N_G(K)=\N_M(K)$ and $\Fix(K)=\{ \a\}$.
\qed

\vskip 3mm
Now we are ready to state the main lemma in this subsection.

\begin{lem}\label{PSL}
Let $G=\PSL(2,q)$ where $q=p^n$ and $G_\a=M=\PSL(2,p^m)$ where either $\frac nm$ is an odd prime or $n=2m$ and $p=2$.  Consider the primitive right multiplication action of $G$ on $\O$. Then
$$\begin{array}{lll}|\G|&=&p^{n-m}\frac{p^{2n}-1}{p^{2m}-1}-1-(p^{n-m}-1)(p^m+1)-\frac 12(\frac{p^n-1}{p^m-1}-1)p^m(p^m+1)\\&-&\frac 12(\frac{p^n+(-1)^{\frac{n}{m}-1}}{p^m+1}-1)p^m(p^m-1).\end{array}$$
In particular, if $n=2m$ and $p=2$, then $|\G|=0$, that is $b(G)>2$.\end{lem}
\demo From the above two lemmas, we may just consider two cases: $K\cong\ZZ_p^m$
 and $K\cong \ZZ_{\frac {p^m\pm 1}d}.$

 Let $K\cong \ZZ_p^m$.  By Lemma~\ref{sl}.(i),  $K$ fixes $p^{n-m}$ points. Since $K$ fixes $|\N_M(K):K|=\frac{p^m-1}d$ points in each suborbit with the point stabilizer $K$, we have $(p^{n-m}-1)/ \frac{p^m-1}d$ such suborbits, while the union of all such suborbits contains $(p^{n-m}-1)/\frac{p^m-1}d \cdot \frac {|M|}{p^m}=(p^{n-m}-1)(p^m+1)$ points of $\O$.

Let $K\cong \ZZ_{\frac {p^m\pm 1}d}$.  By  Lemma~\ref{sl}.(ii),   $K$ fixes $\frac{p^n\pm 1}{p^m\pm 1}$ (resp. $\frac{p^{2m}-1}{p^m\pm 1}$)  if  $\frac nm$ is odd (resp.  $n=2m$) points. Since $K$ fixes $|\N_M(K):K|=2$ points in each suborbit with a point stabilizer $K$, we have $\frac 12(\frac{p^n\pm 1}{p^m\pm 1}-1)$ (resp. $\frac 12(\frac{p^{2m}-1}{p^m\pm 1}-1))$  for  such suborbits, while the union of all such suborbits contains $\frac 12(\frac{p^n\pm 1}{p^m\pm 1}-1)p^m(p^m\mp 1)$  (resp. $\frac 12(\frac{p^{2m}-1}{p^m\pm 1}-1)p^m(p^m\mp 1)$) points of $\O$.

Since each nontrivial subgroup $K_1$ of  $K\cong \ZZ_p^m$ or $\ZZ_{\frac {p^m\pm 1}d}$ fixes the same point set with that of $K$, we do not have suborbits with a point stabilizer $H_\b=K_1.$

Since $|\O|=p^{n-m}\frac{p^{2n}-1}{p^{2m}-1}$, we have
$$\begin{array}{lll}|\G|&=&p^{n-m}\frac{p^{2n}-1}{p^{2m}-1}-1-(p^{n-m}-1)(p^m+1)-\frac 12(\frac{p^n-1}{p^m-1}-1)p^m(p^m+1)\\&-&\frac 12(\frac{p^n+(-1)^{\frac{n}{m}-1}}{p^m+1}-1)p^m(p^m-1).\end{array}$$
One can check that if $n=2m$ and  $p=2$, then $|\G|=0$.  \qed

\subsubsection{$G=\PGL(2,q)$, $M=\PGL(2,p^m)$}
\begin{lem}\label{PGLPGL}
Suppose that  $G=\PGL(2,q)$ and $M=\PGL(2,p^m)$ with $q=p^n=p^{mr}$ for some odd prime $r$. Then for any  $T$-regular suborbit $\Delta$, there exists another $T$-regular suborbit $\Delta'$  such that $\Delta \cup \Delta'$ is a $G$-regular suborbit.
\end{lem}
\demo Let $G=\PGL(2,q)$ and $M=\PGL(2,p^m)$. Set $\O=[G:M]$ and $\a=M$. Let $\G$ be the union of regular  orbits of $M_0$ on $\O $. In what follows, we show that for any regular suborbit $\Delta$ of $T$, there exists another
  regular suborbit $\Delta'$ of $T$  such that $\Delta \cup \Delta'$ is a  regular suborbit of $G$.

   On the contrary, let $\b=\a^g\in \O$ where  $(M_0)_\b=1$, but $K:=M_\b\cong \ZZ_2$ and $K\not\le T$. Then
   $\N_G(K)=\D_{2(p^n\pm 1)}$ and $\N_M(K)=\D_{2(p^m\pm 1)},$  for $p^n\equiv \mp 1(\mod 4).$ Moreover, $G$ has only one conjugacy class of involutions in  $G\setminus T$.
Suppose that $\ZZ_{\frac{p^m\pm1}{2}}\cong H_1\leq T$, $K\nleq H_1$ and $H_1\leq H_2\cong \ZZ_{p^m\pm1}$.
 Since  $\N_G(H_i)=\N_G(K)$ and $\N_M(H_i)=\N_M(K)$, it follows that $H_i$ and $K$ have the same set of fixed points for $i=1,2$.
 This is contrary to $(M_0)_\b=1$.
  \qed

\subsubsection{$G=\PGammaL(2,q)$, $M=\PGL(2,p^m):\ZZ_n$ }\label{PGammaLmr}
Suppose that $r\geq 3$. The main result of this subsection is the following Lemma~\ref{PGammaL}.
\begin{lem}\label{PGammaL}
Let $G=\PGammaL(2,p^{mr})$, $M=\PGL(2,p^m).\ZZ_{mr}$ with $r$ an odd prime, $p^m>2$. Consider the primitive permutation representation of $G$ on $\O=[G:M]$ with
right coset multiplication action. Then the Saxl graph $\Sigma(G)$ has diameter $2$.
\end{lem}

Let $G=\PGammaL(2,p^{mr})$, $\a=M=\PGL(2,p^m).\lg f\rg $ with $n=mr$ and $r$ is an odd prime.
Let $S:=\PGL(2,p^{n})$ and $M_1=\PGL(2,p^m)$.
To show the lemma, we need to know which regular suborbit $\a^{sM_1}$  of $S$  is fixed by a nontrivial subgroup of $\lg vf\rg$ for some $v\in \PGL(2,p^m)$.
Suppose that $\a^{sM_1}$  is such a regular suborbit.  Then it must be fixed by an element $uf'$ of prime order $r'$ with some $u\in \PGL(2,p^m)$, where
$f'=f^{\frac n{r'}}$.  Then
$\a^{s}$ is fixed by a subgroup $K$ of order $r'$. Since $M_1$ acts regularly on this   suborbit, it follows that $K=\lg uf'\rg $ for some $u\in M_1$.

In what follows, let $n'(r')$ be the number of the set of regular suborbits of $S$ which are fixed by $K$. We shall deal with two cases
according to   $r'=r$ and $r'\neq r$, separately.

\begin{lem} \label{r_1=r} Suppose that $r'=r$. Then $n'(r')\le \frac 14(r'-1)^2p^m(p^m+1)/(p^{m}-1)$. \end{lem}
\demo   Suppose that $r'=r$.  Then $r'$ is odd and  $M\geq M_1\times \lg f'\rg $. Now we are dealing with  four cases: $r'\nmid |M_1|$, $r'=p$, $r'\di (p^m-1)$ and $r'\di (p^m+1)$.

\vskip 3mm

(1) Case 1: $r'\nmid |M_1|$.

\vskip 3mm
In this case, every subgroup of the form $\lg uf'\rg $ of order $r'$  in $M$ is conjugate to $\lg f'\rg $.
Since $\N_G(\lg f'\rg )=\N_M(\lg f'\rg)=M_1:\lg f\rg$, the only fixed point of $f'$
is $\a$ and $n'(r')=0.$

\vskip 3mm
(2) Case 2: $r'=p$.

\vskip 3mm
Now $\Gal(\FF_q/\FF_{p^m})=\lg f'\rg $. Since  $p$ is odd and $\sum_{i=0}^{p-1}1^{f^{i}}=0$, there exists $x\in \FF_q$ such that $x-x^{f'}=1$ according to Proposition~\ref{Gal}. Set
$$u_0=\overline{\small{\left(
\begin{array}{cc}
1&1 \\
0&1 \\
\end{array}
\right)}}\quad {\rm and}\quad g=\overline{\small{\left(
\begin{array}{cc}
1&x \\
0&1 \\
\end{array}
\right)}} .$$
Since all the elements of order $p$ in $M_1$ are conjugate to $u_0$ in $M_1$, every subgroup of the form $\lg uf'\rg $ of order $r'$  in $M$ is conjugate to either
$K_0=\lg f'\rg $ or $K_1=\lg u_0f'\rg $.
As $M_1\times \lg f' \rg\leq M$, it follows that $K_0$ and $K_1$ are not conjugate in $M$.
Moreover, one may check that
$(g^{-1})^{f'}g=\overline{\small{\left(
\begin{array}{cc}
1&x-x^{f'} \\
0&1 \\
\end{array}
\right)}}=u_0$, that is, $f'^g=f'u_0$. In other words, all $G$-conjugacy subgroups of $K_0$ (resp. $K_1$)  which are contained in $M$ form two conjugacy classes of $M$, with the respective representatives $K_0$ and $K_1$. Therefore,  by Proposition~\ref{man}, we have
 $$|\Fix ( K_0)|=\frac {|\N_G(K_0)|}{|\N_M(K_0)|}+\frac {|\N_G(K_1)|}{|\N_M( K_1)|}=\frac{|M_1|n}{|M_1|n}+\frac{|M_1|n}{p^mn}=1+\frac{p^m(p^{2m}-1)}{p^m}=p^{2m}.$$
With the same arguments, we have $|\Fix (K_1)|=p^{2m}.$

For $K_0$, we have $n'(r')=0$. Otherwise,   take a regular suborbit $\Delta$  of $S$ which is fixed setwise by  $K_0$. Considering the action of $M_1K_0$ on $\Delta$,  we know that
 $K_0$ fixes $|\N_{M_1K_0}(K_0):K_0|=|M_1|=p^m(p^{2m}-1)$ points, which is impossible, noting that  $K_0$ fixes  $p^{2m}$ points.

As for $K_1$,  in each regular suborbit $\Delta$  of $S$  which is fixed setwise by $K_1$, we get that $K_1$ fixes
 $|\N_{M_1K_1}(K_1):K_1|=p^{m}$ points. Since   $K_1$ fixes $p^{2m}$ points where one of them is $\a$,   the number $n'(r')$ is no more than  $\frac{p^{2m}-1}{p^{m}}$.

\vskip 3mm
(3) Case $3$: $r'\mid p^m-1$.
\vskip 3mm
 Set $p^m=kr'+1$. Then $\frac{p^n-1}{p^m-1}=\frac{(kr'+1)^{r'}-1}{kr'}\equiv 0(\mod r')$.
 Again, we have $\Gal(\FF_q/\FF_{p^m})=\lg f'\rg $.
  For any $a\in \FF_{p^m}$ such that $|a|=r'$,  since for any $1\le i\le r'$,
  $$\prod_{j=0}^{r'-1}(a^i)^{f'^j}=\prod_{j=0}^{r'-1}(a^i)^{p^{mj}}=
  a^{i\frac{p^{n}-1}{p^m-1}}=1,$$
   by Proposition~\ref{Gal}, we may pick up $x\in \FF_q$ such that $xx^{-f'}=a^{if'}$. Set
$$u_0=\overline{\small{\left(
\begin{array}{cc}
1&0 \\
0&a \\
\end{array}
\right)}}\quad {\rm and}\quad g(i)=\overline{\small{\left(
\begin{array}{cc}
1&0 \\
0&x \\
\end{array}
\right)}} .$$

Since all the subgroups  of order $r'$ in $M_1$ are conjugate to $\lg u_0\rg $ in $M_1$, every subgroup of the form $\lg uf'\rg $ of order $r'$  in $M$ is conjugate to either
$K_0=\lg f'\rg $ or $K_i=\lg u_0^if'\rg $, where $1\le i\le \frac12(r'-1)$. Using the same arguments as in Case $2$, we get that $K_0$ and $K_i$ are not conjugate in $M$. Since $u_0^i$ and $u_0^j$ are conjugate in $M$ if and only if $i\equiv-j(\mod r')$, it follows that $K_i$ and $K_j$ are conjugate in $M$ if and only if $i\equiv-j(\mod r')$. Moreover, one may check that $f'^{g(i)}=u_0^if'$, which implies all subgroups of order $r'$ in $M$ are conjugate in $G$.
Now by Proposition~\ref{man}, we have
$$|\Fix (K_0) |=\frac {|\N_G (K_0) |}{|\N_M (K_0) |}+\sum_{i=1}^{\frac 12(r'-1)}\frac {|\N_G (K_i) |}{|\N_M( K_i )|}$$
$$=\frac{|M_1|n}{|M_1|n}+\frac 12(r'-1)\frac{|M_1|n}{(p^m-1)n}=\frac 12(r'-1)p^{m}(p^m+1)+1.$$
With the same arguments,  $|\Fix (K_i) |=\frac 12(r'-1)p^{m}(p^m+1)+1.$

With the same arguments as in Case 2,  $n'(r')=0$ for $K_0$.  As for $K_i$, except for $\a$, $K_i$ fixes $\frac 12(r'-1)p^{m}(p^m+1)$ points.  In each suborbit of $S$ fixed by $K_i$, $K_i$ fixes
 $|\N_{M_1K_i}(K_i):K_i|=(p^m-1)$ points and so $n'(r')\le \frac 14(r'-1)^2p^m(p^m+1)/(p^{m}-1)$.

\vskip 3mm
(4) Case $4$: $r'\mid p^m+1$.
\vskip 3mm
(i) $p=2$: Let $n=mr'$, where $r'$ is an odd prime.  Let $\FF_{2^m}^*=\lg \d\rg $, $\FF_{2^{2m}}=\FF_{2^m}(\zeta)$  where $\zeta^{2^m}=\zeta+1$ and $\zeta^2=\zeta+\th $ for some $\th \in \FF_{2^m}$ (by Proposition~\ref{Gal2}). Since $\zeta\not\in \FF_{2^n}$, $\FF_{2^{2n}}=\FF_{2^n}(\zeta).$
 Take a generator $x+y\zeta$ of $\FF_{2^{2n}}^*$. Then there is an injective homomorphism $\phi$:
   $\FF_{q^2}^*/\FF_q^*\to \PGL(2,q)$, given by $$(x+y\zeta)\FF_q^* \longmapsto s:={{\overline{\left(
\begin{array}{cc}
x&y\th \\
y&x+y\\
\end{array}
\right)  }}}.$$
 Set
$$w=\overline{\left(
\begin{array}{cc}
1& \th\\
0 &1
\end{array}
\right)},\quad  t=\overline{\left(
\begin{array}{cc}
1& 1\\
0 &1
\end{array}
\right)} \quad {\rm and} \quad c=fw.$$ One can check that $\sum_{i=0}^{m-1}\th^{f^i}=\zeta^{f^m}+\zeta=1$, noting $p=2$.
Then
$$s^c=\overline{\left(
\begin{array}{cc}
1& -\th\\
0 &1
\end{array}
\right)}\overline{\left(
\begin{array}{cc}
x^2&y^2\th^2 \\
y^2&x^2+y^2\\
\end{array}
\right)}\overline{\left(
\begin{array}{cc}
1& \th\\
0 &1
\end{array}
\right)}=\overline{\left(
\begin{array}{cc}
x^2+y^2\th&y^2\th \\
y^2&x^2+y^2+y^2\th\\
\end{array}
\right)}
=s^2,$$
$$c^m=(fw)^m=f^m{{ \overline{\left(
\begin{array}{cc}
1& \sum_{i=0}^{m-1}\th^{f^i}\\
0&1 \\
\end{array}
\right)  }}}=f^mt, \quad  c^{2m}=f^{2m}, \quad c^{n}=t \quad {\rm and}\quad c^{2n}=1.$$
Suppose  $r'\di 2^m+1$. Then  $r'\nmid 2^m-1$ and $r'\di\di \frac{2^n+1}{2^m+1}$. Set $r'^j\di\di (2^m+1)$. Then
$r'^{j+1}\di\di (2^n+1)$. Set $s_0=s^{\frac{2^n+1}{r'^{j+1}}}$ and now $R:=\lg s_0\rg :\lg c^{2m}\rg $ is a Sylow $r'-$subgroup of $\lg s \rg:\lg c^{2m}\rg $.
Suppose that $s_0^{c^{2m}}=s_0$.  Then $(2^{2m}-1)\equiv 0(\mod r'^{j+1})$, that is  $(2^m+1)\equiv 0(\mod r'^{j+1})$,
a contradiction. So $R$ is a nonableian $r'-$group which contains a cyclic maximal subgroup.  By Lemma~\ref{maxcyclic} again, all the noncenter subgroups of order $r'$ are conjugate to $\lg  c^{2m}\rg =\lg f'\rg $ in $R$, noting that $f'^2=f^{2m}=c^{2m}$.

\vskip 3mm
(ii) $p$ is odd:  Let $\FF_{p^{m}}^*=\lg \d\rg $ and  $\FF_{p^{2n}}=\FF_{p^n}(\zeta)$ where $\zeta^2=\th :=\d^{\frac{p^m-1}{2^l}}$ and $p^m-1=2^lk$  for an odd integer $k$. In particular, $|\th |=2^l$.   Take a generator $x+y\zeta$ of $\FF_{p^{2n}}^*$. Set
$$s={{ \overline{\left(
\begin{array}{cc}
x&y\th \\
y&x \\
\end{array}
\right)  }}},\quad  t={{ \overline{\left(
\begin{array}{cc}
1&0 \\
0&-1 \\
\end{array}
\right)  }}} \quad {\rm and} \quad w={{ \overline{\left(
\begin{array}{cc}
1&0 \\
0&\th^{\frac{1-p}2} \\
\end{array}
\right)  }}}.$$
Set $c=fw$.
  From the proof  of Lemma~\ref{M}, we have $s^c=s^{p}$, $c^{m}=f^{m}t$, $c^{2m}=f^{2m}$ and $c^n=t$, for any  odd prime divisor $r'$ of $n$.

Suppose  $r'\di p^m+1$ and  set $r'^j\di\di (p^m+1)$. Then
$r'^{j+1}\di\di (p^n+1)$. Set $s_0=s^{\frac {p^n+1}{r'^{j+1}}}$ and again we consider a Sylow $r'$-subgroup $R:=\lg s_0\rg :\lg c^{2m}\rg $ of $\lg s\rg:\lg c^{2m}\rg $.
Then with the similar analyses in (i), we know that $R$ is nonabelian and all the noncenter subgroups of order $r'$  are conjugate to $\lg f'\rg $ in $R$, where $\lg f'\rg=\lg f^{2m}\rg$, as  desired.

\vskip 3mm

Let $\lg u_0\rg$ be a subgroup of order $r'$ of $\lg s\rg $. Since all the subgroups  of order $r'$ in $M_1$ are conjugate to $\lg u_0\rg $ in $M_1$, every subgroup of the form $\lg uf'\rg $ of order $r'$  in $M$ is conjugate to either
$K_0=\lg f'\rg $ or $K_i=\lg u_0^if'\rg $. Moreover, by (i) and (ii), every $K_i$ is conjugate to $\lg f'\rg $ in $G$.
Using the same arguments as in Case $2$ and Case $3$, we know that $K_0$ and $K_i$ are not conjugate to each other in $M$ for $1\leq i\leq \frac{1}{2}(r'-1)$.
Therefore, by Proposition~\ref{man}, we have
$$|\Fix (K_0) |=\frac {|\N_G(K_0)|}{|\N_M(K_0)|}+\sum_{i=1}^{\frac 12(r'-1)}\frac {|\N_G(K_i)|}{|\N_M( K_i )|}$$
$$=\frac{|M_1|n}{|M_1|n}+\frac 12(r'-1)\frac{|M_1|n}{(p^m+1)n}=\frac 12(r'-1)p^{m}(p^m-1)+1.$$
With the same arguments, $|\Fix (K_i) |=\frac 12(r'-1)p^{m}(p^m-1)+1.$

Similarly,  $n'(r')=0$ for $K_0$.  As for $K_i$, except for $\a$, $K_i$ fixes $\frac{1}{2}(r'-1)p^{m}(p^m-1)$ points.  In each suborbit of $S$ fixed by $K_i$, $K_i$ fixes
  $|\N_{M_1K_i}(K_i):K_i|=(p^m+1)$ points and so $n'(r')\le \frac 14(r'-1)^2p^m(p^m-1)/(p^{m}+1)$.

\vskip 3mm
Taking the maximal $n'(r')$ in the above  four cases, we get the result as desired.\qed

\begin{lem} \label{r'nr} Suppose that $r'\ne r$. Then  $n'(r')\le \frac{p^{m_1r}(p^{2m_1r}-1)}{[p^{m_1}(p^{2m_1}-1)]^2},$ where $m=m_1r'$.
\end{lem}
\demo   Suppose that $r'\ne r$. Set $m=m_1r'$. Then $M=M_1:\lg f\rg $. In what follows, we shall prove that  every subgroup of the form $\lg uf'\rg $ of order $r'$  in $M$ is conjugate to $\lg f'\rg $ where $u\in M_1$. Since $\N_G(\lg f'\rg )=\PGL(2,p^{\frac n{r'}})\lg f\rg$ and
$\N_M(\lg f'\rg)=\PGL(2,p^{m_1})\lg f\rg$,  it follows that $\lg f'\rg$  fixes $\frac{p^{m_1r}(p^{2m_1r}-1)}{p^{m_1}(p^{2m_1}-1)}$ points. In each regular  suborbit of $S$ fixed by $\lg f'\rg$, we get that $\lg f'\rg$ fixes
 $|\N_{M_1\lg f'\rg}(\lg f'\rg ):\lg f'\rg |=p^{m_1}(p^{2m_1}-1)$ points. Therefore,
$$n'(r')\le \frac{p^{m_1r}(p^{2m_1r}-1)}{[p^{m_1}(p^{2m_1}-1)]^2}.$$

In fact, the conclusion is clearly  true for  $r'\nmid |M_1|$.  So we are going to deal with the other three cases separately.

\vskip 3mm
(1) Case $1$: $r'=p$.

\vskip 3mm
Note that $\Gal(\FF_{p^m}/\FF_{p^{m_1}})=\lg f'\rg $.   For any $a\in \FF_{p^m}$ with trace $\Tr(a)=0$, we have   $\Tr(a^{f'})=0$.
By Proposition~\ref{Gal},  there exits  $x\in \FF_{p^m}$ such that $x-x^{f'}=a^{f'}$. Set
$$u=\overline{\small{\left(
\begin{array}{cc}
1&a \\
0&1 \\
\end{array}
\right)}}\quad {\rm and}\quad g=\overline{\small{\left(
\begin{array}{cc}
1&x \\
0&1 \\
\end{array}
\right)}} .$$
Up to conjugacy,   we may set $K=\lg uf'\rg $, for some $a\in \FF_{p^m}$. Then $|uf'|=p$ if and only if
$uu^{f'}u^{f'^2}\cdots u^{f'^{p-1}}=1$, that is, $\Tr(a)=a+a^{f'}+\cdots +a^{f'^{r'-1}}=0$.
Then $$f'^g=g^{-1}g^{f'^{-1}}f'=
{{ \overline{\left(
\begin{array}{cc}
1&x^{f'^{-1}}-x\\
0&1 \\
\end{array}
\right)}}}f'=uf'.$$
Therefore, every subgroup of the form $\lg uf'\rg $ of order $p$  in $M$ is conjugate to $\lg f'\rg $.

\vskip 3mm
\vskip 3mm
(2) Case $2$: $r'\mid p^m-1$.
\vskip 3mm
We will discuss two cases where   $r'$ is either odd or even, separately.

Suppose that $r'$ is odd. Since $r'\di p^{m_1r'}-1$,  we have $p^{m_1r'r}-1\equiv 0(\mod r')$ and then $p^{m_1r(r'-1)}p^{m_1r}\equiv1 (\mod r')$.
Therefore, $p^{m_1r}\equiv1 (\mod r')$ and $\frac{p^n-1}{p^{m_1r}-1}\equiv 0(\mod r')$.
Take  $a\in \FF_{p^m}$ such that $|a|=r'$.  Since $\Gal(\FF_{p^m}/\FF_{p^{m_1}})=\lg f'\rg $ and for any $1\le i\le r'$, it follows that
  $$\Norm(a^i)=\prod_{j=0}^{r'-1}(a^i)^{f'^j}=\prod_{j=0}^{r'-1}(a^i)^{p^{m_1rj}}=
  a^{i\frac{p^{n}-1}{p^{m_1r}-1}}=1.$$
By    Proposition~\ref{Gal} again,   there exists $x\in \FF_{p^m}$ such that $xx^{-f'}=a^{if'}$. Set
$$u_0=\overline{\small{\left(
\begin{array}{cc}
1&0 \\
0&a \\
\end{array}
\right)}}\quad {\rm and}\quad g(i)=\overline{\small{\left(
\begin{array}{cc}
1&0 \\
0&x \\
\end{array}
\right)}} .$$
Since all the subgroups  of order $r'$ in $M_1$ are conjugate to $\lg u_0\rg $ in $M_1$, every subgroup of the form $\lg uf'\rg $ of order $r'$  in $M$ is conjugate to $K_i=\lg u_0^if'\rg $, where $0\le i\le \frac12(r'-1)$. Moreover, one may check that $f'^{g(i)}=u_0^if'$, as desired.

Suppose that $r'=2$.
%Note that any Sylow $2$-subgroup of $M_1$ is contained in a group which is isomorphic to $D_{2(p^m-1)}$.
Let $\lg\e \rg$ be a Sylow $2$-subgroup of $\FF_{p^m}^*$ with order $2^l$. Set
$$u_0=\overline{\small{\left(
\begin{array}{cc}
1&0 \\
0&\e \\
\end{array}
\right)}},\quad v=\overline{\small{\left(
\begin{array}{cc}
0&1 \\
1&0 \\
\end{array}
\right)}}\quad {\rm and}\quad y=\overline{\small{\left(
\begin{array}{cc}
1&1 \\
\e^{f'}&\e \\
\end{array}
\right)}}.$$
Consider the 2-group $L=\lg u_0, v, f'\rg$. The involutions in $L\setminus M_1$  have the form either $u_0^jf'$  or $u_0^jvf'$.

Let $d=vf'$, it is easy to check that $u_0^d=u_0^{-p^{\frac{n}{2}}}$ and $f'^{y}=d$. Let $e=\frac{p^{\frac{n}{2}}+1}{(p^{\frac{n}{2}}+1)_2}$ and $e'=\frac{p^{\frac{n}{2}}-1}{(p^{\frac{n}{2}}-1)_2}$, where $(p^{\frac{n}{2}}\pm 1)_2$ denote the largest $2$-power divisor  of $p^{\frac{n}{2}}\pm1$.
Observe that $(u_0^{ei}d)^2=1$ if and only if $ei=(p^{\frac{n}{2}}+1)k$ for some integer $k$ and $(u_0^{e'i'}f')^2=1$ if and only if $e'i'=(p^{\frac{n}{2}}-1)k'$ for some integer $k'$. Note that $(u_0^{-1})^{d^{-1}}=(u_0^{-1})^d=u_0^{p^{\frac{n}{2}}}$ and $u_0^{f'}=u_0^{p^{\frac{n}{2}}}$. Then we may get
$$d^{u_0^{-k}}=u_0^{(p^{\frac{n}{2}}+1)k}d=u_0^{ei}d \quad {\rm and}\quad  f'^{u_0^{k'}}=u_0^{(p^{\frac{n}{2}}-1)k'}f'=u_0^{e'i'}f'.$$
Thus every   involution of $L\setminus M_1$   is conjugate to $f'$.

\vskip 3mm
(3) Case $3$: $r'\mid p^m+1$.
\vskip 3mm
If $r'=2$, then $m=2m_1$ is even and  every Sylow $2$-subgroup of $M_1$  is contained in a subgroup of $M_1$  isomorphic to $\D_{2(p^m-1)}$, which is included in the previous case.
 So suppose that  $r'$ is odd.

\vskip 3mm
(i) $p=2$:  As before, $m=m_1r'$.  Let $\FF_{2^{m_1}}^*=\lg \d\rg $, $\FF_{2^{2m_1}}=\FF_{2^{m_1}}(\zeta)$  where $\zeta^{2^{m_1}}=\zeta+1$ and $\zeta^2=\zeta+\th $ for some $\th \in \FF_{2^{m_1}}$ (by Proposition~\ref{Gal2}). Since $\zeta\not\in \FF_{2^n}$, $\FF_{2^{2m}}=\FF_{2^m}(\zeta)$ and $\FF_{2^{2n}}=\FF_{2^n}(\zeta).$
  As in Lemma~\ref{r_1=r}.(4), take a generator $x+y\zeta$ of $\FF_{2^{2n}}^*$ , set
  $$s:={{\overline{\left(
\begin{array}{cc}
x&y\th \\
y&x+y\\
\end{array}
\right)  }}}, \quad
 w=\overline{\left(
\begin{array}{cc}
1& \th\\
0 &1
\end{array}
\right)},\quad  t=\overline{\left(
\begin{array}{cc}
1& 1\\
0 &1
\end{array}
\right)} \quad {\rm and} \quad c=fw.$$
 Then we have $s^c=s^2$ and   $\sum_{i=0}^{m_1-1}\th^{f^i}=\a^{f^{m_1}}+\a=1$. Moreover,  $$\sum_{i=0}^{m_1r-1}\th^{f^i}=\sum_{j=0}^{r-1}(\sum_{i=0}^{m_1-1}\th^{f^i})^{f^{jm_1}}=r\cdot 1=1.$$
Then
$$c^{m_1r}=(fw)^{m_1r}=f^{m_1r}{{ \overline{\left(
\begin{array}{cc}
1& \sum_{i=0}^{m_1r-1}\th^{f^i}\\
0&1 \\
\end{array}
\right)  }}}=f^{m_1r}t, ~\, c^{2m_1r}=f^{2m_1r},~  c^{n}=t~{\rm and} ~c^{2n}=1.$$
Suppose  $r'\di 2^m+1$. Then $2^m+1\equiv 2^{m_1}+1\equiv 0(\mod r')$, that is $r'\di 2^{m_1}+1$.  Set $r'^j\di\di (2^{m_1}+1)$. Then
$r'^{j+1}\di\di (2^m+1)$.  Set $s_0= s^{\frac{2^n+1}{r'^{j+1}}}$ and   $R:=\lg s_0\rg :\lg c^{2m_1r}\rg $ is a Sylow $r'-$subgroup of $M_1:\lg c^{2m_1r}\rg $.
Suppose that  $[s_0, c^{2m_1r}]=1$.
Then $(2^{2m_1r}-1)\equiv 0(\mod r'^{j+1})$, that is  $(2^{m_1r}+1)\equiv 0(\mod r'^{j+1})$. However, set $2^{m_1}=r'^jk-1$ where $(k,r')=1$.
Then $2^{m_1r}+1=(r'^jk-1)^r+1=r'^jk'$ for some $k'$ coprime to $r'$. So  $[s_0, c^{2m_1r}]\ne 1$ and then
$R$ is a nonableian $r'-$group which contains a cyclic maximal subgroup.  By Lemma~\ref{maxcyclic} again, all the noncenter subgroups of order $r'$ are conjugate to $\lg  c^{2m_1r}\rg =\lg f'^2\rg=\lg f'\rg $ in $R$, noting that $f'^2=f^{2m_1r}=c^{2m_1r}$.

\vskip 3mm
(ii) $p$ is odd:  Let $\FF_{p^{m_1}}^*=\lg \d\rg $ and  $\FF_{p^{2m_1}}=\FF_{p^{m_1}}(\zeta)$ with $\zeta^2=\th :=\d^{\frac{p^{m_1}-1}{2^l}}$, where $p^{m_1}-1=2^lk$  for an odd integer $k$. In particular, $|\th |=2^l$.   Take a generator $x+y\zeta$ of $\FF_{p^{2n}}^*$. Set
$$s={{ \overline{\left(
\begin{array}{cc}
x&y\th \\
y&x \\
\end{array}
\right)  }}},\quad  t={{ \overline{\left(
\begin{array}{cc}
1&0 \\
0&-1 \\
\end{array}
\right)  }}} \quad {\rm and} \quad w={{ \overline{\left(
\begin{array}{cc}
1&0 \\
0&\th^{\frac{1-p}2} \\
\end{array}
\right)  }}}.$$
Set $c=fw$. As in  Lemma~\ref{M}, $s^c=s^{p}$. Moreover, we get that

$$c^{m_1r}=f^{m_1r}{{ \overline{\left(
\begin{array}{cc}
1&0 \\
0&(\th^{\frac{1-p}2})^{\frac{1-p^{m_1r}}{1-p}} \\
\end{array}
\right)  }}} =f^{m_1r}{{ \overline{\left(
\begin{array}{cc}
1&0 \\
0&-1 \\
\end{array}
\right)  }}}=f^{m_1r}t.$$
So $c^{2m_1r}=f^{2m_1r}$ and $c^n=t$ for any  odd prime divisor $r'$ of $n$.
The remaining proof is completely the same as in (i), just replacing $p=2$ by odd $p$.
\qed

\vskip 3mm
Now we are ready to prove Lemma~\ref{PGammaL}.
\vskip 3mm \f {\bf Proof of Lemma~\ref{PGammaL}} In Lemmas~\ref{PSL} and \ref{PGLPGL}, we already determine the number of regular suborbits of $T=\PSL(2,q)$ and $S=\PGL(2,q)$, relative to $\a=\PGL(2,p^m):\lg f\rg $
on $\O$. In order to determine the regular orbits  of $M$, we first get a bound for  the number
 $n'$, which is the size of the set of regular suborbits of $S$ which are fixed by a nontrivial subgroup of $M\setminus S$ with prime order.

If $m\ge 2$, we set $m=r_1^{e_1}r_2^{e_2}\cdots r_l^{e_l}$, where $2\le r_1<r_2 \cdots <r_l$ are  primes and $e_i\ne 0$. Then $m\ge 3^{l-1}r_1$, that is $l\le 1+\log_3\frac m{r_1}.$
  Then
$$n'\le n'(r)+ \sum_{i=1}^ln'(r_i) \, {\rm for}\,  r\nmid m \quad {\rm and}\quad n'\le n'(r)+ \sum_{r_i\neq r}n'(r_i) \, {\rm for}\,  r\mid m.$$
Particularly, we have $n'=n'(r)$ for $m=1$. Moreover, noting $p^m\ge 3$,  by Lemmas~\ref{r_1=r} and ~\ref{r'nr}, we have
$$n'(r)\leq \frac 14(r-1)^2p^m(p^m+1)/(p^{m}-1)\leq \frac 12(r-1)^2p^m\quad {\rm and}\quad$$
 $$ n'(r_i)\le \frac{p^{m_ir}(p^{2m_ir}-1)}{[p^{m_i}(p^{2m_i}-1)]^2}\le\frac{16}{9}p^{3m_ir-6m_i}\le \frac{16}{9}p^{3m_1r-6m_1},$$ where $m=m_ir_i$.
 Since $n=rr_1^{e_1}r_2^{e_2}\cdots r_l^{e_l}$, we get  that

$$ \begin{cases}
n'\le \frac 12(r-1)^2p^m,  & \text{if } m=1 \\
n'\le \frac 12(r-1)^2p^m+l \frac{16}{9}p^{3m_1r-6m_1}\le \frac 12(r-1)^2p^m+(1+\log_3m_1)\frac{16}{9}p^{3m_1r-6m_1}, & \text{if } m> 1
\end{cases}
 .$$

Let $t_1,t$ be the  number of regular suborbits of $S$ and $G$, respectively. Note that $t_1|M_1|$ has been computed in the proof of Lemma \ref{PSL}.  Then $t|M|\ge (t_1-n')|M_1|$. Now
$$\begin{array}{lll} &&t|M|-\frac{|\O|}{2}\ge t_1|M_1|-n'|M_1|-\frac{|\O|}2\\ &\geq&\frac{1}{2p^m(p^{2m}-1)}(p^{3n}+2p^{n+4m}-2p^{n+3m}+2p^{n+m}+p^n-2p^{2m}-2p^{3m}-2p^{5m}+2p^{4m})\\
&&-n'p^m(p^{2m}-1)\\
&=&\frac{1}{2p^m(p^{2m}-1)}p^{3n}(A+(p^{n+4m}-2p^{n+3m})p^{-3n}+(p^{n+4m}+2p^{n+m}+p^n+2p^{4m}-2p^{2m}\\
&& -2p^{3m}-2p^{5m})p^{-3n})\ge \frac{1}{2p^m(p^{2m}-1)}p^{3n}A>0,  \\
  \end{array}$$
where $A:=1-2n'p^{2m}(p^{2m}-1)^2p^{-3n}\ge 1-2n'p^{6m-3n}=1-A_1-A_2>0$, with
 $$A_1=2p^{6m-3n}\frac 12(r-1)^2p^m=\frac{(r-1)^2}{p^{m(3r-7)}}\le \frac{(r-1)^2}{3^{3r-7}}<\frac 12;$$
  and either $m=1$ with $A_2=0$, or   $m\ge 2$  with
  $$A_2=2p^{6m-3n}(1+\log_3m_1)\frac{16}{9}p^{3m_1r-6m_1}\le \frac{32}{9}\frac{(1+\log_3m_1)}{p^{3(r-2)(m-m_1)}}
  \le \frac{32}{9}\frac{(1+\log_3m_1)}{p^{3m_1}}<\frac 12.$$
Therefore,  $t|M|>\frac{|\O|}{2}$ and so  the Saxl graph of $\PGammaL(2,q)$ has diameter 2.
  \qed

\subsection{$M_0=\D_{\frac{2(q-1)}{d}}$ where $d=(2,q-1)$}
If $p=2$, then $b(G)=2$ with $\soc(G)=\PSL(2,q)$ if and only if $G=\PGL(2,q)$ by \cite[Lemma 4.7]{B}. Followed by Example $2.5$ in \cite{BG}, we have that the Saxl graph $\Sigma(\PGL(2,q))$ has diameter $2$. Hence, from now on, we assume that $p$ is an odd prime and  $q\neq 5,7,9,11$ in this subsection.
Also by \cite[Example 2.5]{BG}, we have that the Johnson graph $\J(|\O|,2)$ is a subgraph of both $\Sigma(T)$ and $\Sigma(\PGL(2,q))$ for $M_0=\D_{\frac{2(q-1)}{d}}$. And so in these cases, the Saxl graph has diameter $2$.

Let $G=\PSigmaL(2,p^{n})$, $M=\D_{q-1}.\lg f\rg $ where $f$ is induced by  the field automorphism  $\sigma: a\mapsto a^p$ for any $a\in \FF_{p^{n}}^*$ and $q=p^n$, where $n\ge 1$.
Let $T:=\PSL(2,p^{n})$ and $M_0=\D_{q-1}$.
Set $\O'=[T:M_0]$ and $\O=[G:M]$.  Consider   the primitive permutation representation of $G$ on $\O$. Then the action of $T$ on $\O$ is
equivalent to that on $\O'$.

Now $\O$ may be identified with the set of $2$-subsets of $\PG(1,q)$ and $G$ induces a natural action on $\O$. The projective line $\PG(1,q)$ can be identified with $\FF_q\cup \{\infty\}$ by mapping $\lg (x,y)\rg $ to
$\frac{y}{x}$. Let $\a:=\{0,\infty \}$
 and $\FF_{p^{n}}^*=\lg \th \rg $. Set

$$s={{\small{ \overline{\left(
\begin{array}{cc}
1&0 \\
0&\th^2 \\
\end{array}
\right)  }}}}\quad {\rm and}\quad v={{{\small \overline{\left(
\begin{array}{cc}
0&-1 \\
1&0 \\
\end{array}
\right)  }}}}.$$
Then $M:=G_{\a }=\lg  s,v,f \rg=\lg s,v\rg:\lg f \rg \cong \D_{q-1}.\ZZ_n$ and  $M_0:=T_\a =\lg s,v \rg$.

Set $G_1:=\lg T,\delta f^i\rg$ and $G_2:=T:\lg f^i \rg$ where $\frac{n}{(n,i)}$ is even and $\delta={{\small{ \overline{\left(
\begin{array}{cc}
1&0 \\
0&\th \\
\end{array}
\right)  }}}}$ . In what follows, for $L\in \{ T, G, G_1, G_2\}$,  by  $X_L(\a)$ and $Y_L(\a)$ we denote the union of all irregular and regular suborbits of $L$ relative to $\a$, respectively. Set $S^*=(\FF_q^*)^2$  and $N=\FF_q^*\setminus S^*.$
 %  Note that Johnson graph is a subgraph of $\Sigma(T) $ whose  vertices $\{u,w\}$  ($\ne \a$)  are fixed by $M\setminus M_0$, provides  one of   $\{u,w\}$ is either 0 or $\infty$,  for $\Sigma(\PGL(2,q))$ only have one  regular suborbit.

\begin{lem}\label{fixed}
Let $u,w\in\FF_q$ and $u\neq w$. Then $\{u, w \}\in X_T(\a)$ if and only if $\frac{w}{u}=-\th^{2i}$ for some $\th^i\in\FF_q^*$.
\end{lem}
\demo
First, suppose $\{u,w\}\in X_T(\a)$.
Then $\{ u,w\}^{s^i}=\{ u, w\}$  or  $\{u, w\}^{s^iv}=\{u,w\}$ for some $i\in \ZZ_{q-1}$,
that is, $\{u,w\}=\{\th^{2i}u, \th^{2i}w\}$ or $\{-\th^{-2i}u^{-1}, -\th^{-2i}w^{-1}\}$.
 Then we get either $w=-u$ for $q\equiv 1(\mod 4)$ or $\frac{w}{u}=-\th^{-2i}u^{-2}$.
 In both cases, $\frac{w}{u}=-\th^{2i'}$ for some $\th^{i'}\in\FF_q^*$.

 Conversely, suppose $w=-u\th^{2i}$ for some $\th^i\in\FF_q^*$ and $u=\th^{k}$ for some integer $k$.
 Then $\{ u ,w\}^{s^{-i-k}v}=\{-u^{-1}\th^{2i}u^2,-w^{-1}\th^{2i}u^2\}=\{w,u\}$, which means $\{u,w\}\in X_T(\a)$.
 \qed

\begin{lem}\label{delta}
Let $u,w\in\FF_q$ and $u\neq w$. Then $\{u, w \}\in Y_T(\a)\setminus Y_G(\a)$ if and only if  either   one of $u$ and $w$ is  $0$ or $\infty$;  or   $u, w\in \FF_q^*$,
 and $\frac{w}{u}\in \FF_{p^r}^*$, where $r\di n$ and $r<n$.
   \end{lem}
\demo
Suppose $\{u,w\}\in Y_T(\a)$.
First, no loss,   assume  $u$   is $0$ or $\infty$. Then    $\{0,w\}^{f^rs^i}=\{0,w\}$ or  $\{\infty,w\}$ where $\th^{2i}=w^{1-p^r}$for some integer $i$.
So  in what follows, let $u, w\not\in \a$.

Set  $\{u,w\}^g=\{u,w\}$ for some $g\in M\setminus M_0$ where $\frac wu\in -N$.
Set $g=f^ks^iv^j$ where $1\leq k < n$ and $j\in \{0,1\}$.
Then $\{u,w\}^g=\{\th^{2i}u^{p^{k}},\th^{2i}w^{p^k}\}=\{u,w\}$ for $j=0$ and $\{u,w\}^g=\{-\th^{-2i}u^{-p^{k}},-\th^{-2i}w^{-p^k}\}=\{u,w\}$ for $j=1$.
By solving it, we get
$$\begin{cases}
\th^{\epsilon2i}u^{\epsilon p^{k}}=\epsilon u  \\
\th^{\epsilon2i}w^{\epsilon p^k}=\epsilon w
\end{cases}
\quad {\rm or} \quad
\begin{cases}
\th^{\epsilon2i}u^{\epsilon p^{k}}=\epsilon w  \\
\th^{\epsilon2i}w^{\epsilon p^k}=\epsilon u
\end{cases}$$
where $\epsilon=\pm1$. Then $(\frac{w}{u})^{p^{k}+1}=1$ or $(\frac{w}{u})^{p^{k}}=\frac{w}{u}$.

 To show 
$(\frac{w}{u})^{p^{k}}=\frac{w}{u}$, equivalently $(\frac{w}{u})^{p^r}=\frac{w}{u}$ where $r=(n,k)$,  we  prove   $(\frac{w}{u})^{p^k+1}\neq 1$ for any $1\leq k<n$. On the contrary, assume that  $(\frac{w}{u})^{p^k+1}=1.$

(1) Suppose that  $(n,2k)=(n,k)$. From    $(p^k+1, p^k-1)=2$ and $(p^n-1,p^{2k}-1)=p^{(n,2k)}-1=p^{(n,k)}-1$, we get
 $(\frac 12(p^n-1),\frac12(p^{k}-1)(p^k+1))=\frac12 (p^{(n,k)}-1) (p^{(n,k)}+1),$  and then
  $(p^n-1,p^k+1)=2$. Then $(\frac{w}{u})^{p^k+1}=1$ if and only if $(\frac{w}{u})^2=1$, that is $\frac wu=\pm 1$, a contradiction.

  (2) Suppose that $(n,2k)=2(n,k)$. Then   we may set $k=(n,k)k_1$ and $n=(n,k)n_1$, where $k_1$ is odd, $(n_1,k_1)=1$, $n_1$ is even, as $(n,2k)=2(n,k)$.
   Since
      $$p^{nk_1}-1=((p^k+1)-1)^{n_1}-1=(\sum_{i=0}^{n_1}\binom{n_1}{i} (p^k+1)^{i}(-1)^{n_1-i})-1\equiv  0(\mod 2(p^k+1))$$
     and $k_1$ is odd, we have $p^n-1\equiv 0(\mod 2(p^k+1)_2)$ where $(p^k+1)_2$ denotes the largest $2$-power of $p^k+1$. This contradicts to  $\frac wu\in -N=N$ and   $(\frac wu)^{p^k+1}=1$.

Conversely, suppose $(\frac{w}{u})^{p^{r}}=\frac{w}{u}$ and $u,w\in \FF_q^*$. Then $w^{p^{r}-1}=u^{p^{r}-1}=\th^{2i}$ for some $\th^i\in\FF_q^*$. Then
$$\{u,w\}^{f^rs^{-i}}=\{\th^{-2i}u^{p^r},\th^{-2i}w^{p^r}\}=\{\th^{-2i}u\th^{2i},\th^{-2i}w\th^{2i}\}=\{u,w\},$$
that is $\{u, w \}\in Y_T(\a)\setminus Y_G(\a)$. This completes the proof of this lemma.
\qed

\vskip 3mm
From what we have discussed above, we get that the Johnson graph is a subgraph of the Saxl graph $\Sigma(T)$ but not a subgraph of $\Sigma(\lg T,f^j  \rg)$ for any $1\leq j<n$.

\vskip 3mm
\begin{lem}\label{Dsigma}
Consider the primitive permutation representation of $G$ on $\Omega$ with right coset multiplication action. Then the Saxl graph $\Sigma(G)$ has diameter $2$.
\end{lem}
\demo
The lemma  has  been confirmed by Magma for  $q\leq 79$ and so  we assume that $q>79$. Set $S^*=(\FF_q^*)^2$  and $N=\FF_q^*\setminus S^*.$
From the previous lemma, we know
  $$\begin{array}{ll} X_G(\a):=&\{\{u,w\}\mid \frac{w}{u}=-s~or~\frac{w}{u}=-t\in\FF_{p^r}^*~{\rm where}~s\in S^*, r\di n, r\not\equiv0(\mod n)\\ &~{\rm and}~w,u\in \FF_q^* \}
                    \cup \{\{\infty, w\}\mid w\in \FF_q^*\} \cup \{\{u, 0\}\mid u\in \FF_q^*\}.\end{array}$$
%and $\widehat{Y(\a)}:=\{\{u,w\}\mid \frac{w}{u}\neq-s~and~\frac{w}{u}\neq-t\in\FF_{p^k}^*~where~s\in S^*~and~k\not\equiv0(\mod n) \}$.

We shall deal with three cases, separately.

\vskip 3mm
(1) $\b_0=\{u_0,w_0\}\in X_G(\a)$ where $\{0,\infty\} \cap \{u_0,w_0\}=\emptyset$.

 Set  $g={{\small{ \overline{\left(
\begin{array}{cc}
1&u_0 \\
\lambda&\lambda w_0 \\
\end{array}
\right)  }}}}$ such that  $\{0,\infty\}^g=\{u_0,w_0\}$, where $\lambda:=\frac{1}{w_0-u_0}$ and $w_0=-b_0u_0$ with $b_0\in S^*\cup\FF_{p^r}^*$. Then
$X_G(\b_0)=X_G(\a)^g=A\cup B\cup C$, where
 $$\begin{array}{l}A=\{\frac{u_0+\lambda uw_0}{1+u\lambda},\frac{u_0+\lambda w_0w}{1+\lambda w}\}\mid
 \frac{w}{u}=-s~or~\frac{w}{u}=-t\in\FF_{p^r}^*~
 ~s\in S^*, r\di n, r\not\equiv0(\mod n)\}\\
 B=\{\{w_0,\frac{u_0+
 \lambda w_0w}{1+\lambda w}\}\mid w\in \FF_q^*\},\quad  C=\{\{\frac{u_0+\lambda uw_0}{1+u\lambda},u_0\}\mid u\in \FF_q^*\}.
  \end{array}$$
 We try to find a pair $\{\frac{u_0+\lambda u'w_0}{1+\lambda u'},\frac{u_0+\lambda w_0w'}{1+\lambda w'}\}\notin A\cup B\cup C$ for some $\{u',w'\}$ where $w'=-x u'$ for a primitive element $x\in \FF_q^*$ such that
  \begin{eqnarray}\label{eta-1}
 -\th=\frac{u_0+\lambda w_0w'}{1+\lambda w'}\frac{1+u'\lambda}{u_0+\lambda u'w_0}.
\end{eqnarray}
Then the pair $\{\frac{u_0+\lambda u'w_0}{1+\lambda u'},\frac{u_0+\lambda w_0w'}{1+\lambda w'}\}\in Y_G(\a)\setminus X_G(\b_0)$ and so it is adjacent to both $\a$ and $\b_0$.

Now, Eq(\ref{eta-1}) is equivalent to
 \begin{eqnarray}
\label{x-1}
F_{x}(u'):=(b_0x+\th b_0x)u'^2-(b'+b_0xb'-\th b'b_0-\th xb')u'+(b'^2+\th b'^2)=0.
\end{eqnarray}
Set
\begin{eqnarray}\label{delta(x)}
\Delta(x):=b'^2[(\th-b_0)^2x^2-2(b_0\th^2+\th+\th b_0^2+4\th b_0+b_0)x+(b_0\th-1)^2]
\end{eqnarray}
where  $b':=u_0(1+b_0)$ and $b_0\neq \th$ because $b_0\in S^*\cup \FF_{p^r}^*$.
For given $b_0$, $u_0$ and some $x\neq 0$,  Eq(\ref{x-1}) has solutions for $u'$  if and only if $\Delta(x)\in S^*$.
By the corollary in \cite[page 173]{DW},  there exists a primitive element $\eta$ such that $\Delta(\eta)\in S^*$ if $q\not\in \{121, 169\}$,
 noting that $q$ is not a prime. Moreover, for $q\in \{121, 169\}$, Magma shows that such primitive element $\eta$ still exists  if $b_0\in \FF_q^*\setminus\{-1\}$.
Now take $x=\eta$  so that    $u'=\frac{(b'+b_0\eta b'-\th b'b_0-\th \eta b')\pm (\Delta(\eta ))^{\frac{1}{2}} }{2(b_0\eta +\th b_0\eta )}$ and $w'=-\eta u'$.
Then we shall show the pair $\{\frac{u_0+\lambda u'w_0}{1+\lambda u'},\frac{u_0+\lambda w_0w'}{1+\lambda w'}\}\not\in A\cup B\cup  C.$

Suppose that   $\{\frac{u_0+\lambda u'w_0}{1+\lambda u'},\frac{u_0+\lambda w_0w'}{1+\lambda w'}\}\in A.$ Then    $\{\frac{u_0+\lambda u'w_0}{1+\lambda u'},\frac{u_0+\lambda w_0w'}{1+\lambda w'}\}= \{\frac{u_0+\lambda uw_0}{1+\lambda u},\frac{u_0+\lambda w_0w}{1+\lambda w}\}$ for some $u, w=-bu$ where
$b\in S^*\cup\FF_{p^r}^*$. Since $\frac{u_0+\lambda w_0w'}{1+\lambda w'}\frac{1+u'\lambda}{u_0+\lambda u'w_0}=-\th =\frac{u_0+\lambda w_0w}{1+\lambda w}\frac{1+u\lambda}{u_0+\lambda uw_0}$, we get
  $\frac{u_0+\lambda u'w_0}{1+\lambda u'}=\frac{u_0+\lambda uw_0}{1+\lambda u}$ and   $\frac{u_0+\lambda w_0w'}{1+\lambda w'}=\frac{u_0+\lambda w_0w}{1+\lambda w}.$ These two equations imply  that
    $u'=u$ and $w'=w$, which in turn gives $\eta=b$, a contradiction. Therefore, $\{\frac{u_0+\lambda u'w_0}{1+\lambda u'},\frac{u_0+\lambda w_0w'}{1+\lambda w'}\}\not\in A.$

Suppose that   $\{\frac{u_0+\lambda u'w_0}{1+\lambda u'},\frac{u_0+\lambda w_0w'}{1+\lambda w'}\}\in B.$ Then  either   $\frac{u_0+\lambda u'w_0}{1+\lambda u'}=w_0$ or $\frac{u_0+\lambda w_0w'}{1+\lambda w'}=w_0$. Both cases give $w_0=u_0$, a contradiction. So  $\{\frac{u_0+\lambda u'w_0}{1+\lambda u'},\frac{u_0+\lambda w_0w'}{1+\lambda w'}\}\not\in B.$
Similarly, one has  $\{\frac{u_0+\lambda u'w_0}{1+\lambda u'},\frac{u_0+\lambda w_0w'}{1+\lambda w'}\}\not\in C.$

In summary, $\{\frac{u_0+\lambda u'w_0}{1+\lambda u'},\frac{u_0+\lambda w_0w'}{1+\lambda w'}\}\not\in A\cup B\cup  C,$ as desired.

 \vskip 3mm
  (2) $\b_0=\{0,w_0\}\in X_G(\a)$ for any $w_0\in \FF_q^*$.

Set $g={{\small{ \overline{\left(
\begin{array}{cc}
1&0 \\
\lambda&\lambda w_0 \\
\end{array}
\right)  }}}}$ for $\lambda=\frac{1}{w_0}$ such that  $\{0,\infty\}^{g}=\{0,w_0\}$. Then
$X_G(\b_0)=X_G(\a)^g=A\cup B\cup C$, where
 $$\begin{array}{l}A=\{\{\frac{\lambda uw_0}{1+u\lambda},\frac{\lambda w_0w}{1+\lambda w}\}\mid \frac{w}{u}=-s~{\rm or}~\frac{w}{u}=-t\in\FF_{p^r}^*~
 {\rm where}~s\in S^*, r\di n, r\not\equiv0(\mod n)\},\\
 B=\{\{w_0,\frac{\lambda w_0w}{1+\lambda w}\}\mid w\in \FF_q^*\},\quad  C=\{\{\frac{\lambda uw_0}{1+u\lambda},0\}\mid u\in \FF_q^*\}.
  \end{array}$$
 Note that $\frac{\lambda uw_0}{1+u\lambda}$ or $\frac{\lambda w_0w}{1+\lambda w}$ is $\{\infty\}$ if $1+u\lambda=0$ or $1+\lambda w=0$, respectively.

We will prove $\{\frac{1}{\th}\frac{(\th-\th^p)w_0}{\th^p+1},-\frac{(\th-\th^p)w_0}{\th^p+1}\}\notin A\cup B\cup C$, that is,  $\{\frac{1}{\th}\frac{(\th-\th^p)w_0}{\th^p+1},-\frac{(\th-\th^p)w_0}{\th^p+1}\}\in Y_G(\a)\setminus X_G(\b_0)$.

Suppose $\{\frac{1}{\th}\frac{(\th-\th^p)w_0}{\th^p+1},-\frac{(\th-\th^p)w_0}{\th^p+1}\}\in A$. Then $\{\frac{1}{\th}\frac{(\th-\th^p)w_0}{\th^p+1},-\frac{(\th-\th^p)w_0}{\th^p+1}\}=\{\frac{\lambda uw_0}{1+u\lambda},\frac{\lambda w_0w}{1+\lambda w}\}$ for some $u, w=-bu$ where $b\in S^*\cup \FF_{p^r}^*$. Therefore, no loss, let $\frac{\lambda w_0w}{1+\lambda w}\frac{1+u\lambda}{\lambda uw_0}=-\th$. Then $(w_0b+\th w_0b)u^2+(bw_0^2-\th w_0^2)u=0$ and $u=\frac{(\th-b)w_0}{b(\th+1)}$ because $u\neq 0$.
So $\{\frac{\lambda uw_0}{1+u\lambda},\frac{\lambda w_0w}{1+\lambda w}\}=\{\frac{1}{\th}\frac{(\th-b)w_0}{b+1},-\frac{(\th-b)w_0}{b+1}\}$
and $\frac{1}{\th}\frac{(\th-\th^p)w_0}{\th^p+1}=\frac{1}{\th}\frac{(\th-b)w_0}{b+1}$. Therefore, $b=\th^p$, a contradiction.

Suppose $\{\frac{1}{\th}\frac{(\th-\th^p)w_0}{\th^p+1},-\frac{(\th-\th^p)w_0}{\th^p+1}\}\in B$. Then either $\frac{1}{\th}\frac{(\th-\th^p)w_0}{\th^p+1}=w_0$ or
$-\frac{(\th-\th^p)w_0}{\th^p+1}=w_0$. Both cases give $\th=-1$, a contradiction. And clearly $\{\frac{1}{\th}\frac{(\th-\th^p)w_0}{\th^p+1},-\frac{(\th-\th^p)w_0}{\th^p+1}\}\notin C$ because $0\notin \{\frac{1}{\th}\frac{(\th-\th^p)w_0}{\th^p+1},-\frac{(\th-\th^p)w_0}{\th^p+1}\}$.

In summary, $\{\frac{1}{\th}\frac{(\th-\th^p)w_0}{\th^p+1},-\frac{(\th-\th^p)w_0}{\th^p+1}\}\notin A\cup B\cup C$, as desired.

\vskip 3mm
(3) $\b_0=\{u_0, \infty \}\in X_G(\a)$ for any $u_0\in \FF_q^*$.

Because $\infty^{v}=0$, $0^{v}=\infty$ and $v$ is an automorphism of the Saxl graph $\Sigma(G)$, now we can also say that  $\a=\{0,\infty\}$ and $\b_0=\{u_0, \infty\}$ have a common neighbor in the Saxl graph $\Sigma(G)$ for any $u_0\in \FF_q^*$.
 \qed

\begin{lem}
Consider the primitive permutation representation of $G_1:=\lg \PSL(2,q), \delta f^i \rg$ on $\Omega'$ with right coset multiplication action where $\frac{n}{(n,i)}$ is even. Then the Saxl graph $\Sigma(G_1)$ has diameter
$2$.
\end{lem}
\demo
Set $G_1:=\lg T,\delta f^i\rg$ and $G_2:=T:\lg f^i \rg$ where $\frac{n}{(n,i)}$ is even and $\delta={{\small{ \overline{\left(
\begin{array}{cc}
1&0 \\
0&\th \\
\end{array}
\right)  }}}}$ . Then $q\equiv 1(\mod 4)$ and $\{u,w\}\in Y_{T}(\a)$ if and only if $\frac{w}{u}=\th^{2l+1}$ for some $\th^l\in \FF_q^*$ by Lemma \ref{fixed}.
By a similar analysis to Lemma \ref{delta} and Lemma \ref{fixed}, we get $$ Y_{G_2}(\a)=\{\{u,w\}\mid \frac{w}{u}=\th^{2l+1},\quad (\frac{w}{u})^{p^{ik}}\neq \frac{w}{u}\quad {\rm and} \quad ik\not\equiv 0(\mod n) \}.$$
Obviously, $Y_{G_2}(\a)\subset Y_{T}(\a)$. Now we will prove $Y_{G_2}(\a)\subset Y_{G_1}(\a)$. On the contrary, we assume that $\{u,w\}^{\delta^k f^{ik}s^j}=\{u,w\}$ or $\{u,w\}^{\delta^{k} f^{ik}s^jv}=\{u,w\}$ where $ik\not\equiv
0 (\mod n)$ for some $\{u,w\}\in Y_{G_2}(\a)$.  Then we get that
$$\{u,w\}=\{\th^{kp^{ik}+2j}u^{p^{ik}},\th^{kp^{ik}+2j}w^{p^{ik}}\}\quad {\rm or}\quad  \{u,w\}=\{-\th^{-(kp^{ik}+2j)}u^{-p^{ik}},-\th^{-(kp^{ik}+2j)}w^{-p^{ik}}\}.$$
By solving these two equations, we have
$$
(\frac{w}{u})^{p^{ik}+1}=1
\quad {\rm or} \quad
(\frac{w}{u})^{p^{ik}}=\frac{w}{u}
.$$
Since $\{u,w\}\in Y_{G_2}(\a)$, it follows $(\frac{w}{u})^{p^{ik}}\neq\frac{w}{u}$.
And by what we have analyzed in Lemma \ref{delta}, we have $(\frac{w}{u})^{p^{ik}+1}\neq 1$.
So $Y_{G_2}(\a)\subset Y_{G_1}(\a)$. That is, the Saxl graph $\Sigma(G_2)$ is a subgraph of $\Sigma(G_1)$. Followed by Lemma \ref{Dsigma}, we get that the Saxl graph $\Sigma(G_1)$ has diameter $2$.
\qed

 \vskip 3mm
The following lemma follows from the previous lemmas and Lemma $4.7$ in \cite{B}. By \cite[Lemma 4.7]{B}, we have that $b(G)\leq 3$ for $\soc(G)=\PSL(2,q)$ and $M_0=\D_{\frac{2(q-1)}{d}}$. Moreover, the equality holds if and only if $\PGL(2,q)<G$.

\begin{lem}\label{d(q-1)}
Let $\soc(G)=\PSL(2,q)$ and $M_0=\D_{\frac{2(q-1)}{d}}$, where $q\geq 5$ and $q\ne 5,7, 9,11$. Let $f$ be the field automorphism of $\FF_{p^n}$ and $\PGL(2,q)=\lg \PSL(2,q),\delta \rg$. Then $b(G)=d(\si)=2$ if and only if one of the following occurs:

\begin{enumerate}
\item[\rm(i)] $G=\PSL(2,2^n)$;
\item[\rm(ii)] For odd $q$, $G\in \{\PSL(2,q),\PGL(2,q),\PSL(2,q):\lg \d f^{i} \rg,\PSL(2,q):\lg f^{j} \rg\}$ with $\frac{n}{(n,i)}$ is even and $j$ is an integer.
\end{enumerate}

\end{lem}

\subsection{$M_0\in\{A_4, S_4, A_5,\ZZ_p^n:\ZZ_{\frac{q-1}{d}},\PGL(2,p^{\frac{n}{2}})\}$}\label{some}
Now we have five cases for $M_0$. We shall solve the first three cases by  probabilistic methods and  show that $b(G)>2$ for the other two cases.

\begin{lem}\label{GLD} Suppose that  $\soc(G)=\PSL(2,q)$ and  $M_0=A_4$. Then $q=p$ and $b(G)=2$ if and only if $q\ge 11$. Moreover, the Saxl graph has diameter $2$ for $b(G)=2$.
 \end{lem}
  \demo Since $M_0=A_4$, $q=p\equiv \pm 3 (\mod 8)$. Moreover, if $q=p\not\equiv \pm 1(\mod 5)$, $M_0$ is maximal in $T$; otherwise,    $M_0\le A_5\le T$.   In both cases, it suffices to deal with the case: $G=\PGL(2,q)$ and
  $M=S_4$.

 Now $\mathcal{P}^*(M)=\{\lg x_1\rg, \lg x_2\rg, \lg x_3\rg\}$, with order $2$, $2$ and $3$, respectively.
  Since $\N_M(\lg x_1\rg)\cong \D_8$,  $\N_M(\lg x_2\rg)\cong \D_4$ and $\N_M(\lg x_3\rg)\cong \D_6$, and $\lg x_1\rg $ and $\lg x_2\rg $ are not conjugate in $\PGL(2,q)$, we have
   $$
   \begin{array}{lll}\hat{Q}(G)&=&(\frac{(|x_1|-1)}{|N_M(\lg x_1\rg)|}\times \frac{|N_G(\lg x_1\rg)|}{|N_M(\lg x_1\rg)|}+\frac{(|x_2|-1)}{|N_M(\lg x_2\rg)|}\times\frac{|N_G(\lg x_2\rg )|}{|N_M(\lg x_2\rg)|}+\frac{(|x_3|-1)}{|N_M(\lg x_3\rg)|}\times\frac{|N_G(\lg x_3\rg )|}{|N_M(\lg x_3\rg)|})\frac{|M|}{|\O|}\\ \\
&&\leq(\frac{1}{8}\times\frac{2(q+1)}{8}+\frac{1}{4}\times\frac{2(q+1)}{4}+\frac{2}{6}\times\frac{2(q+1)}{6})
  \times\frac{24^2}{q(q^2-1)}=\frac{154}{q(q-1)}<\frac{1}{2},
 \end{array} $$
 for all $q\geq 19$. By Proposition~\ref{smallp},   $d(\si)$=2 if $q\ge 19$.
   Moreover, if $q=5$, then $b(G)=1$; if $q=11, 13, 19$, then by Magma we know that $b(G)=2$ and $d(\si)=2$.
 \qed

\begin{lem}\label{S_4}
 Suppose that  $\soc(G)=\PSL(2,q)$ and  $M_0=S_4$. Then $q=p\equiv \pm 1(\mod 8)$ and $b(G)=2$ if and only if $q\ge 17$. Moreover,  the Saxl graph has diameter $2$ for $b(G)=2$.
 \end{lem}
\demo In this case, $G=\PSL(2,q)$ nad $M=S_4$. Now  $\mathcal{P}^*(M)=\{\lg x_1\rg, \lg x_2\rg, \lg x_3\rg\}$ with the respective order 2, 2 and 3, and $\N_M(\lg x_1\rg)\cong \D_8$,  $\N_M(\lg x_2\rg )\cong \D_4$ and $\N_M(\lg x_3\rg )\cong \D_6$.
 Moreover, $\lg x_1\rg $ and $\lg x_2\rg $ may be  conjugate in $G$ or may not, depending on $q$. So for $i=1,2$, $|\Fix(\lg x_i\rg)|\le \frac{q+1}{8}+\frac{q+1}{4}=\frac{3(q+1)}8$.
 Therefore,
 $$\hat{Q}(G)\leq(\frac{1}{8}\times\frac{3(q+1)}{8}+\frac{1}{4}\times\frac{3(q+1)}{8}+\frac{2}{6}\times\frac{q+1}{6})
 \times\frac{24^2\times2}{q(q^2-1)}=\frac{226}{q(q-1)}<\frac{1}{2}$$  provided $q\geq 23$.
  If  $q=7,9$, then  $b(G)\neq 2$; and  if $q=17$, then $d(\si)=2$ by Magma.
\qed

\begin{lem}\label{A_5}
Suppose  that $\soc(G)=\PSL(2,q)$ and  $M_0=A_5$. Then  either
\begin{enumerate}
\item[\rm(i)] $q=p\equiv \pm 1(\mod 5)$  and $G=\PSL(2,q)$; or
  \item[\rm(ii)] $q=p^2\equiv -1(\mod 5)$ for  an odd prime $p$ and $G=\PSL(2,q)$ or $\PSigmaL(2,q).$
 \end{enumerate}   Moreover,   $b(G)=2$ if and only if  $q\ge 29$, and  the Saxl graph $d(\si)=2 $ for $b(G)=2$.
\end{lem}
\demo (1) Suppose that  $q=p\equiv \pm 1(\mod 5)$, then $G=\PSL(2,q)$.
Now $\mathcal{P}^*(M)=\{\lg x_1\rg, \lg x_2\rg, \lg x_3\rg\}$, with the respective order 2, 3, and 5,  $\N_M(\lg x_1\rg)\cong \D_4$,  $\N_M(\lg x_2\rg )\cong \D_6$ and $\N_M(\lg x_3\rg )\cong \D_{10}$. Therefore,
  $$\hat{Q}(G)\leq(\frac{1}{4}\times\frac{q+1}{4}+\frac{2}{6}\times\frac{q+1}{6}+\frac{4}{10}\times\frac{q+1}{10})\times\frac{60^2\times 2}{q(q^2-1)}
  =\frac{1138}{q(q-1)}<\frac{1}{2}$$  provided $q\ge 49$.
  As for  $q=9,11,16,19$, $|\O|<60$ and so $b(G)\neq 2$.   Checking by Magma,  we know that if $q=29$ (resp. 31), then $G$ has 1 (resp. 2) regular suborbit(s);  if $q=41$, then  $|\Gamma|=420>\frac{1}{2}|\O|$ and therefore,   $b(G)=d(\Sigma)=2$ in these cases.

\vskip 3mm
(2) Suppose $q=p^l\equiv \pm 1(\mod 5)$ for $l\ge 2$. If  $p\equiv \pm 1(\mod 5)$, then $A_5\le \PSL(2,p)<\PSL(2,q)$, which is impossible. So  $p\equiv \pm 2(\mod 5)$ which implies
$p^2\equiv -1(\mod 5)$ and so $l=2$. In this case,  $f$ can fix a subgroup  $A_5$ but $\t$ cannot, where $\PGL(2,q)=\PSL(2,q):\lg \t \rg$. So $G=\PSL(2,q)$ or $\PSigmaL(2,q)$.  It suffices to deal with  $G=\PSigmaL(2,p^2)$ and  $M=S_5$.

 In this case, $\mathcal{P}^*(M)=\{\lg x_1\rg, \lg x_2\rg, \lg x_3\rg, \lg x_4\rg\}$, with the respective order $2$, $2$, $3$ and $5$, and $\N_M(\lg x_1\rg)\cong \D_8$,  $\N_M(\lg x_2\rg )\cong \D_{12}$, $\N_M(\lg x_3\rg )\cong \D_{12}$, $\N_M(\lg x_4\rg )\cong \ZZ_5.\ZZ_4.$ Moreover,
 $$|\Fix(x_i)|\leq \frac{2(p^2-1)}8+\frac{2p(p^2-1)}{12}=\frac {2p^3+3p^2-2p-3}{12}$$ for $i=1,2$.
 Therefore,
  $$\hat{Q}(G)\leq((\frac 1{8}+\frac{1}{12})\times \frac {2p^3+3p^2-2p-3}{12}+
   \frac{2}{12}\times\frac{2(p^2+1)}{12}
   +\frac{4}{20}\times\frac{2(p^2+1)}{20})\times\frac{120^2}{p^2(p^4-1)}
<\frac{1}{2},$$  provided $p\ge 13$.  So $d(\si)=2$ in this case. 
 Checking by Magma, we know that if  $p=3$, then $b(G)\neq 2$;  and if $p=7$,  then $d(\si)=2$.
\qed

\begin{lem} Suppose that either $M_0=\ZZ_p^n:\ZZ_{\frac{q-1}{d}}$ or $M_0=\PGL(2,p^{\frac{n}{2}})$ where $n$ is even and $p$ is an odd prime. Then $b(G)\ne 2$.\end{lem}
\demo
Consider the primitive permutation representation of $T=\PSL(2,q)$ on $\O':=[T:M_0]$.
 Suppose that $M_0=\ZZ_p^n:\ZZ_{\frac{q-1}{d}}$. Then   the action of  $T$ on $\O$  is 2-transitive and  $T_\a\cap T_\b\cong \ZZ_{\frac{q-1}{d}}$, which implies    $b(T)\neq 2$ and $b(G)\ne 2$.

 Suppose $M_0=\PGL(2,p^{\frac{n}{2}})$ be a maximal subgroup of $T$ where $p$ is an odd prime. Consider the primitive action of $T$ on $\O'$. When $q$ is a prime,  it was proved  in \cite[Lemma 4.1]{PX} that   there exist no regular suborbits.  In fact, checking  \cite[Lemma 4.1]{PX}, we find the proof is true  for any prime power $q$. Hence  $b(T)\neq 2$ and  $b(G)\ne 2$.
 \qed
\vskip 3mm
\begin{center}{\large\bf Acknowledgements}\end{center}
\vskip 2mm
We would like to thank the referees for their helpful comments and thank Prof. Keqin Feng who proved  Lemma \ref{feng}.
The first author thanks the supports of Hubei Provincial Natural Science Foundation of China (2020CFB265).
The second author thanks the    supports of the National Natural Science Foundation of China (11671276,12071312,11971248).

\end{document}